\documentclass[11pt,a4paper]{article}
\usepackage{amssymb,amsmath, amsfonts,graphicx,pgfplots,multirow,amsthm}
\usepackage{color}
\usepackage{graphicx}
\usepackage{subcaption}
\usepackage[graphicx]{realboxes}
\usepackage[flushleft]{threeparttable}
\usepackage{lscape}
\usepackage{booktabs}
\usepackage{hyperref}
\usepackage{breqn}
\hypersetup{colorlinks=true, citecolor=blue}
\usepackage[round]{natbib}
\usepackage[margin=1in,footskip=0.5in]{geometry}
\numberwithin{equation}{section}

\newcommand{\BR}{\textbf{R}~}
\newcommand{\br}{\textbf{r}~}
\newcommand{\BS}{\textbf{S}~}
\newcommand{\bs}{\textbf{s}~}

\newcommand{\rsk}{$R_{s;k}$~}
\newcommand{\data}{\textup{data}}

\begin{document}
		\title{\bf Multicomponent stress strength reliability estimation for Pareto distribution based on upper record values}
		\author{Qazi Azhad Jamal$^{a,b}$, Mohd. Arshad$^{a}$ and Nancy Khandelwal$^{a}$\\
		$^{a}$Department of Statistics and Operations Research, Aligarh Muslim University, Aligarh\\
		$^{b}$Department of Mathematics and Statistics, Banasthali Vidyapith, Rajasthan
	}	
	
	\date{}
	\maketitle
	\hrule
\begin{abstract}
In this article, inferences about the  multicomponent stress strength reliability are drawn under the assumption that strength and stress follow  independent Pareto distribution with different shapes $(\alpha_1,\alpha_2)$ and common scale parameter $\theta$. The maximum likelihood estimator, Bayes estimator under squared error and Linear exponential loss function, of multicomponent stress-strength reliability are constructed with corresponding highest posterior density interval for unknown $\theta.$  For known $\theta,$ uniformly minimum variance unbiased estimator and asymptotic distribution of multicomponent stress-strength reliability with asymptotic confidence interval is discussed. Also, various Bootstrap confidence intervals are constructed. A simulation study is conducted to numerically compare the performances of various estimators of multicomponent stress-strength reliability. Finally, a real life example is presented to show the applications of derived results in real life scenarios. 

\vspace{0.5cm}
\end{abstract}\hrule
\section{Introduction}
The stress-strength reliability of a system defines the probability that the system will function properly until the strength exceeds stress. Due to the manufacturing variability and uncertain factors, the strength of the system varies and also when the system is put to use, it is subjected to the stress which is again random in nature. These manufacturing variables and uncertain factors can be used material, production style, humidity, temperature of the environment etc. The genesis of this problem can be seen in \cite{birnbaum1956use}. Later, \cite{birnbaum1958distribution} studied statistical properties of this model.
   Although the model is very simple in nature but its largely applicable in fields of reliability, engineering, manufacturing etc. Since its emergence, various researchers have produced research article on different statistical distributions based on this model. For recent references see, 
    \cite{gunasekera2015generalized}, \cite{wang2018inference}, \cite{sharma2018bayesian}, \cite{ccetinkaya2019stress} and \cite{bai2019reliability}.\\
\indent \cite{bhattacharyya1974estimation} observed that, in several practical scenarios,  the performance of a system depends on more than one component and these components have their own strengths. For example, an aircraft generally contains more than one engines $(k)$ and assume that for takeoff at least $s~ (1\leq s \leq k)$ engines are needed. So, the aircraft will takeoff smoothly, if $s$-out-of-$k$ engines work;
 in engineering,  a power system powering a manufacturing unit has $k$ fuse cutouts arranged in a parallel way. The power system will keep powering the manufacturing unit as long as at least  $s~ (1\leq s \leq k)$ fuse cutouts are working etc. \\ 
\indent Suppose a system, with $k$ identical and independent components, i.e., $X_1,X_2,\ldots,X_k$ works if at least $s(1\leq s\leq k) $ components simultaneously operate. In this framework, the system has common stress $Y,$ and follow distribution function $F.$ The strengths of $s$ components, are random in nature, and follow distribution function $G.$ Then the system reliability $R_{s;k}$ indicates the probability that system does not fail and is given as (see \cite{bhattacharyya1974estimation})
\begin{align}\label{MSS}
	R_{s;k}=\textup{P}(\textup{atleast}~ s~ \textup{of}~ (X_1,X_2,\ldots,X_k) ~\textup{exceeds} ~Y) 
	=\sum_{p=s}^{k}{k \choose p}\int_{-\infty}^{\infty}(1-G(y))^pG(y)^{k-p}dF(y)
\end{align}
 \cite{bhattacharyya1974estimation} were the first to formulate this multicomponent stress-strength (MSS) model. The authors assumed that the random variables associated with the strength and stress of the system are independently exponentially distributed.
  Later, \cite{draper1978bayesian} discussed the Bayesian perspective of the MSS reliability under the assumption that strength and stress follow independent exponential distribution.
   After that \cite{Pandey1991}  discussed the MSS reliability under the assumption that multiple strength and stress follow independent Burr distributions. The authors derived the maximum likelihood estimator (MLE) and Bayes estimator of the system reliability. Also numerical computation suggested that, for small sample, Bayes estimator is better than MLE. 
    \cite{Kizilaslan2015} addressed the MSS reliability of the system for the situation when strength and a common stress follow Weibull distribution. The authors derived MLE  and asymptotic confidence interval of the MSS reliability. Due to complex nature of the MSS reliability, the authors have obtained approximate Bayes estimators by using Lindley's approximation and Markov Chain Monte Carlo (MCMC) method. 
 \cite{kizilaslan2017classical} considered  proportional reverse hazard rate (PRHR) model as the distribution of strength and stress of the system for the classical and Bayesian estimation of the MSS reliability. The author assumed that scale parameter of the distribution $ \lambda$ is common whereas shape parameters are distinct as $\alpha_1$ and $\alpha_2.$ The author derived MLE, asymptotic distribution and asymptotic confidence interval of MSS reliability when all involved parameters are unknown and derived uniformly minimum variance unbiased estimator of reliability when $ \lambda$ is known. The author also obtained Bayes estimator in closed form using independent gamma priors and approximate Bayes estimator using Lindley's approximation and MCMC technique. 
 \cite{kohansal2017estimation} considered the MSS reliability estimation using  progressively Type II censored data under the assumption that strength and stress have Kumaraswamy distribution. The authors considered classical and Bayesian estimation of the system reliability under multiple strength components.
 \cite{ali2018reliability} considered the problem of estimation of reliability of $s$-out-of-$k$ system for non identically distributed strength and stress i.e, strength follows Weibull distribution and stress follows Burr-III distribution. The authors performed classical and Bayesian estimation for estimation of reliability of such system. Further, it is shown that Bayes estimator performs better than MLE.
 \cite{Kzlaslan2018} considered the bivariate Kumaraswamy distribution for the estimation of MSS reliability under classical and Bayesian scheme.  For more references readers  are advised to see 
 \cite{rao2013estimation} ,
 \cite{Nadar2016}, \cite{gadde2017reliability}, 
  and \cite{chaturvedi2018estimation}.\\
 Record values and the associated inferences are of great importance in various practical fields such as hydrology, sports, medical, life testing etc. Record values are defined as the successive extremes in the sequence of random variables.  \cite{chandler1952distribution} originated the concept of record values and later \cite{foster1954distribution} discussed  hypothesis tests for knowing the distribution of record values based on the sum and difference of the upper and lower records in the series. The mathematical formulation of record values is defined as ``Let $X_1$,$X_2$,$\dots$,$X_n$,$\ldots$ be a sequence of independent identically distributed random variables. Define $$U_1 = 1~~ \textup{and}~~ U_{n+1} = \text{min}\{j:j>U_n, X_j>X_{U_n}\}.$$ The sequence $\{ X_{U_n}, n=1,2,\ldots \}$ and $\{U_n,~n=1,2,\ldots\}$ are called upper record values and record times respectively." Let $\textbf{R}=(R_1,R_2,\ldots,R_n)$ be the first upper records from a sequence of random variable having $f(x)$ and $F(x)$ as probability density function $(pdf)$ and distribution function $(df)$ respectively. Then the joint density of $\textbf{R}$ is
\begin{equation}\label{jointrecords}
	f(\textbf{r})=\prod_{i=1}^{n-1}\left(\dfrac{f(r_i)}{1-F(r_i)}\right)f(r_n),
\end{equation}
where $\textbf{r}=(r_1,r_2,\ldots,r_n)$ is the observed value of $\textbf{R}$.
  For a sufficient reading of record values see the books: \cite{ahsanullah1995record} and \cite{arnold1998records} and for their applications, recent research papers:
  \cite{belaghi2015construction}, \cite{kizilaslan2016estimation}, \cite{singh2017bayesian}, \cite{arshad2018interval}, 
    \cite{arshad2018estimation}, \cite{azhad2018TLExp}, \cite{Arshad2019} and \cite{arshad2019kumaraswamy}.
       Statistical research on stress-strength reliability based on record values is growing and various authors have derived results based on such problems. For example, \cite{baklizi2008likelihood,baklizi2014bayesian} considered the estimation of stress-strength reliability for generalized exponential and exponential distribution using record values respectively. \cite{khan2016umvu} considered PRHR model for reliability and the stress-strength reliability estimation for lower record values.
      Recently,  \cite{Rasethuntsa2018} considered the non-identical multicomponent strength with common stress for family of Kumaraswamy distribution using record values. The authors have obtained MLE and UMVUE of the reliability and constructed asymptotic and Bootstrap confidence intervals under classical estimation scheme. The authors have also obtained closed form of Bayes estimators under conjugate priors and approximate Bayes estimators using Lindley approximation and MCMC techniques.
      Recently, \cite{Stat8242} considered the problem of stress-strength reliability estimation of PD based on upper record values. The authors obtained MLE and approximate Bayes estimators for reliability quantity. For more references of reliability estimation of stress-strength under record values framework see
\cite{zakerzadeh2015inference} and \cite{chiang2018inference}.\\
The Pareto distribution (PD) was first introduced by Vilfredo Pareto, who defined it as the distribution of income. The \textit{pdf} of the PD with shape parameter $\alpha$ and scale parameter $\theta$ i.e., $P(\alpha,\theta)$ is
 \begin{equation}\label{pdf pareto}
 f(x;\alpha,\theta)=\alpha \theta^\alpha x^{-(\alpha+1)}, \qquad x\geq \theta,~ \theta >0,~\alpha>0.
 \end{equation}
 The applications of PD can be seen in the field of economics, finance, environmental studies etc. For a detailed literature on PD readers can see the book \cite{arnold2015pareto}. Various statistical inferences of PD are explored by several researchers. For example,
   \cite{han2017bayesian} discussed the E-Bayesian and hierarchical Bayesian of the PD under squared error loss function and asymmetric loss functions.  Further, a monte carlo study is carried out and performances of estimators are considered using mean square errors and bias values. \cite{tripathi2017estimation} considered the improved estimation of shape parameter of PD for unknown scale. Classes of improved estimators using different techniques are also proposed.   For more references on PD see 
   \cite{shafay2017bayesian} and \cite{jiang2018using}.\\
 This article deals with estimation of the MSS reliability defined in \eqref{MSS}, when the underline distribution of the strength and stress follows non identical PD and also the considered data is of record type. Let $X_1,X_2,\ldots X_k$  be the strength of the system which is independently and identically distributed as $P(\alpha_1,\theta)$ and $Y$ be the common stress of the system, distributed as $P(\alpha_2,\theta).$ From \eqref{MSS} and \eqref{pdf pareto}, the MSS reliability is 
 \begin{align}\label{MSSpareto}
 	R_{s;k}&=	\sum_{p=s}^{k}{k \choose p}\int_{-\infty}^{\infty}(1-G(y))^pG(y)^{k-p}dF(y) \nonumber\\
 	&=\sum_{p=s}^{k}{k \choose p}\int_{\theta}^{\infty}\left(\dfrac{\theta}{y}\right)^{p\alpha_1}\left(1-\left(\dfrac{\theta}{y}\right)^{\alpha_1}\right)^{k-p}\dfrac{\alpha_2\theta^{\alpha_2}}{y^{\alpha_2+1}}dy\nonumber\\
 	&=\sum_{p=s}^{k}\sum_{u=0}^{k-p}{k \choose p}{k-p \choose u}(-1)^u\alpha_2\theta^{\alpha_1(p+u)+\alpha_2}\int_{\theta}^{\infty}\frac{1}{y^{\alpha_1(p+u)+\alpha_2+1}}dy\nonumber\\
 	&=\sum_{p=s}^{k}\sum_{u=0}^{k-p}(-1)^u{k \choose p}{k-p \choose u}\left(\dfrac{\alpha_2}{\alpha_1(p+u)+\alpha_2}\right).
 \end{align}
\section{Estimation of $R_{s;k}$ for unknown $\theta$}
In this section, various estimators of \rsk under classical and Bayesian schemes are obtained. 
\subsection{Maximum Likelihood Estimation of \rsk}
Let $\textbf{R}=(R_1,R_2,\ldots,R_n)$ be the first \textit{n} upper record values from $ P(\alpha_1,\theta)$ and  $\textbf{S}=(S_1,S_2,\ldots,S_m)$ be the \textit{m} first upper record values from $P(\alpha_2,\theta)$ independent of \BR. From \eqref{MSS}, \eqref{jointrecords}, and \eqref{pdf pareto}, the log likelihood function of $\alpha_1$, $\alpha_2$ and $\theta$ is 
\begin{multline}\label{loglike}
 	\ln L(\alpha_1,\alpha_2,\theta|\br,\bs)=n\ln\alpha_1+m\ln\alpha_2+(\alpha_1+\alpha_2)\ln\theta-\alpha_1\ln r_n-\alpha_2\ln s_m-\sum_{i=1}^{n}\ln r_i-\sum_{i=1}^{m}\ln s_i\\
 	\alpha_1>0,\alpha_2>0,\theta<\textup{min}\{r_1,s_1\}
 \end{multline}
 where \br and \bs are realizations of \BR and \BS respectively. From \eqref{loglike}, it is easy to derive the MLEs of unknown parameters i.e.,  $\alpha_1,~\alpha_2$ and $\theta$ as
\begin{align*}
\hat{\theta}&=\textup{min}\{r_1,s_1\}\\
\hat{\alpha_1}&=\dfrac{n}{\ln r_n-\ln ({\textup{min}\{r_1,s_1\}})}\\
\hat{\alpha_2}&=\dfrac{m}{\ln s_m-\ln ({\textup{min}\{r_1,s_1\}})}
 \end{align*} 
From the invariance property of MLE, we get the MLE of \rsk as
\begin{equation}\label{mlersk}
	\hat{R}_{s;k}= \sum_{p=s}^{k}\sum_{u=0}^{k-p}(-1)^u{k \choose p}{k-p \choose u}\left(\dfrac{\hat{\alpha_2}}{\hat{\alpha_1}(p+u)+\hat{\alpha_2}}\right).
\end{equation} 
\subsection{Bayesian Inference of \rsk}
In this section, the authors have derived the Bayes estimators of \rsk using Lindley approximation and MCMC simulation technique because of the complex nature of \rsk. For parameters $\alpha_1$, $\alpha_2$ and $\theta$, it is assumed that the prior distributions are independent and follow two parameter gamma distribution i.e,
\begin{equation*}\label{prioralpha}
	\pi(\alpha_i)=\dfrac{b_i^{a_i}}{\Gamma(a_i)}\alpha_i^{a_i-1}e^{-b_i\alpha_i},\quad a_i,b_i,\alpha_i>0,~i=1,2,
\end{equation*} 
 and
 \begin{equation*}\label{theta}
	\pi(\theta)=\dfrac{b_{3}^{a_3}}{\Gamma(a_3)}\theta^{a_3-1}e^{-b_3\theta},\quad a_3,b_3,\theta>0. 	
 \end{equation*}
Due to independent nature of prior, the joint prior distribution is
\begin{equation}\label{jointprior}
	\pi(\alpha_1,\alpha_2,\theta)=\pi(\alpha_1)\pi(\alpha_2)\pi(\theta),\quad \quad	\alpha_1>0,\alpha_2>0,\theta>0.
\end{equation}
This section considers squared error loss (SEL) function  as symmetric and Linear Exponential (LINEX) loss   function as asymmetric loss function for conducting Bayesian study under  this reliability setup. The SEL function is defined as
\begin{equation*}
L_{s}(\delta, \boldsymbol{\lambda})=	(\delta- \boldsymbol{\lambda})^2,~~ \boldsymbol{\lambda}>0,
\end{equation*}
with posterior mean as Bayes estimator $\left(\delta_{s}\right).$ 
 Also, the LINEX loss function is defined as
 \begin{equation*}\label{linexloss}
 	L_{l}(\delta, \boldsymbol{\lambda})=e^{c(\delta- \boldsymbol{\lambda})}-c(\delta- \boldsymbol{\lambda})-1,\qquad c\ne 0
 \end{equation*}
 with corresponding Bayes estimator as
 \begin{equation*}\label{linexbayes}
 	\delta_{l}=-\dfrac{1}{c}\ln \left(E(e^{-c \boldsymbol{\lambda}}|data)\right).
 \end{equation*}
We see that, \rsk has a complex mathematical form and close form of  Bayes estimators for such mathematical functions is a tedious task.
  To resolve this issue, two numerical techniques Lindley approximation(\cite{lindley1980approximate}) and MCMC are employed.  The Bayes estimators of \rsk for SEL and LINEX loss function are defined respectively as
  \begin{align}
  \delta_{s}&=\dfrac{\int\limits_S    R_{s;k}	\pi(\alpha_1,\alpha_2,\theta|\data)d\alpha_1d\alpha_2d\theta}{\int\limits_S 	\pi(\alpha_1,\alpha_2,\theta|\data)d\alpha_1d\alpha_2d\theta}\\ 
  \delta_{l}&=-\dfrac{1}{c}\ln \left(\dfrac{\int\limits_S \textup{exp}\left(-cR_{s;k}\right)	\pi(\alpha_1,\alpha_2,\theta|\data)d\alpha_1d\alpha_2d\theta}{\int\limits_S 	\pi(\alpha_1,\alpha_2,\theta|\data)d\alpha_1d\alpha_2d\theta}\right).
  \end{align}
  where $\int\limits_S$ denotes the triple integral over the set $S=\{ (0,\infty)^2\times(0,\textup{min}(r_1,s_1)) \}.$
\subsubsection{Lindley Approximation}\label{Lindleymethod}
\cite{lindley1980approximate} approximated the ratio of the two integrals using Taylor series expansion, which is fairly applicable in the calculation of the expectation of the posterior densities. The quantity of Bayes estimator required to approximate is
\begin{equation*}
	E(w( \boldsymbol{\lambda})|data)=\dfrac{\int w( \boldsymbol{\lambda})e^{\mathbb{L}( \boldsymbol{\lambda})+\rho( \boldsymbol{\lambda})}d \boldsymbol{\lambda}}{\int e^{\mathbb{L}( \boldsymbol{\lambda})+\rho( \boldsymbol{\lambda})}d \boldsymbol{\lambda}}
\end{equation*}
where $\mathbb{L}$ is the logarithm of the likelihood function, $\rho$ is the logarithm of the prior distribution of $ \boldsymbol{\lambda}$ and $ \boldsymbol{\lambda}$ is the vector of $m$ parameters $ (\lambda_1, \lambda_2,\ldots, \lambda_m).$ Using the Lindley approximation method, $E(w( \boldsymbol{\lambda})|data)$ can be approximated as
 \begin{equation}
	E(w( \boldsymbol{\lambda})|data)\approx w+\dfrac{1}{2}\sum_{i}^m\sum_{j}^m\left(w_{ij}+2w_i\rho_j\right)\sigma_{ij}+\dfrac{1}{2}\left(\sum_{i}^m\sum_{j}^m\sum_{k}^m\sum_{p}^m\mathbb{L}_{ijk}\sigma_{ij}\sigma_{kp}w_p\right)\bigg\arrowvert_{ \boldsymbol{\lambda}=\hat{ \boldsymbol{\lambda}}}
\end{equation}
where $\hat{ \boldsymbol{\lambda}}$ is the MLE of  $ \boldsymbol{\lambda} = (\lambda_1, \lambda_2,\ldots, \lambda_m)$,  \textit{m} denotes the number of unknown parameters, $\mathbb{L}_{ijk}=\partial^3\mathbb{L}/\partial \lambda_i\partial \lambda_j\partial \lambda_k$, $\rho_j=\partial\rho/\partial \lambda_j$ and $w_{ij}=\partial^2w/\partial\lambda_i \partial \lambda_j$ and $\sigma_{ij}$ is the $(i, j)th$ element of the inverse of  matrix $[-\mathbb{L}_{ij}].$ All of the used quantities in this method are obtained at the MLE of $ \boldsymbol{\lambda}$ i.e., $\hat{ \boldsymbol{\lambda}}.$
In our problem, our quantity of interest, \rsk contains three unknown parameters as $\alpha_1$, $\alpha_2$, and $\theta$ i.e, $ \boldsymbol{\lambda} =(\alpha_1,\alpha_2,\theta).$ Therefore, the Lindley approximation is given by
 \begin{equation}\label{lindley}
 	E(w(\boldsymbol{\lambda})|data) \approx w+\left(\sum_{i=1}^{3}w_id_i+d_4+d_5\right)+\dfrac{1}{2}\left(A\sum_{i=1}^{3}w_i\sigma_{1i}+B\sum_{i=1}^{3}w_i\sigma_{2i}+C\sum_{i=1}^{3}w_i\sigma_{3i}\right)
 \end{equation}
 where
 \begin{align*}
 			d_i=&\rho_1\sigma_{i1}+\rho_2\sigma_{i2}+\rho_3\sigma_{i3},\qquad i=1,2,3, \\d_4=&w_{12}\sigma_{12}+w_{13}\sigma_{13}+w_{23}\sigma_{23},\\d_5=&\dfrac{1}{2}
 	\left(w_{11}\sigma_{11}+w_{22}\sigma_{22}+w_{33}\sigma_{33}\right),\\
 	A=&\mathbb{L}_{111}\sigma_{11}+2\mathbb{L}_{121}\sigma_{12}+2\mathbb{L}_{131}\sigma_{13}+2\mathbb{L}_{231}\sigma_{23}+\mathbb{L}_{221}\sigma_{22}+\mathbb{L}_{331}\sigma_{33},\\
 	B=&\mathbb{L}_{112}\sigma_{11}+2\mathbb{L}_{122}\sigma_{12}+2\mathbb{L}_{132}\sigma_{13}+2\mathbb{L}_{232}\sigma_{23}+\mathbb{L}_{222}\sigma_{22}+\mathbb{L}_{332}\sigma_{33},\\
 	C=&\mathbb{L}_{113}\sigma_{11}+2\mathbb{L}_{123}\sigma_{12}+2\mathbb{L}_{133}\sigma_{13}+2\mathbb{L}_{233}\sigma_{23}+\mathbb{L}_{223}\sigma_{22}+\mathbb{L}_{333}\sigma_{33}.
 \end{align*}
Under SEL function, we have
\begin{equation}
w(\alpha_1,\alpha_2,\theta)=\sum_{p=s}^{k}\sum_{u=0}^{k-p}(-1)^u{k \choose p}{k-p \choose u}\left(\dfrac{\alpha_2}{\alpha_1(p+u)+\alpha_2}\right).
\end{equation}
The required quantities for Lindley approximation are     
\begin{align*}
		w_1&=\dfrac{\partial w}{\partial \alpha_1}=\sum_{p=s}^{k}\sum_{u=0}^{k-p}(-1)^{u+1}{k \choose p}{k-p \choose u}\dfrac{\alpha_2(p+u)}{(\alpha_1(p+u)+\alpha_2)^2},\\
	w_2&=\dfrac{\partial w}{\partial \alpha_2}=\sum_{p=s}^{k}\sum_{u=0}^{k-p}(-1)^{u}{k \choose p}{k-p \choose u}\dfrac{\alpha_1(p+u)}{(\alpha_1(p+u)+\alpha_2)^2},\\
	w_{11}&=\dfrac{\partial^2w}{\partial\alpha_1^2}=\sum_{p=s}^{k}\sum_{u=0}^{k-p}(-1)^{u}{k \choose p}{k-p \choose u}\dfrac{2\alpha_2(p+u)^2}{(\alpha_1(p+u)+\alpha_2)^3} \\
		w_{12}&=\dfrac{\partial^2w}{\partial\alpha_1\partial\alpha_2}=\sum_{p=s}^{k}\sum_{u=0}^{k-p}(-1)^{u}{k \choose p}{k-p \choose u}\dfrac{(p+u)(\alpha_2-\alpha_1(p+u))}{(\alpha_1(p+u)+\alpha_2)^3}=w_{21}.\\
		w_{22}&=\dfrac{\partial^2w}{\partial\alpha_2^2}=\sum_{p=s}^{k}\sum_{u=0}^{k-p}(-1)^{u+1}{k \choose p}{k-p \choose u}\dfrac{2\alpha_1(p+u)}{(\alpha_1(p+u)+\alpha_2)^3}\\
		w_3&=\dfrac{\partial w}{\partial \theta}=0, w_{13}=\dfrac{\partial^2 w}{\partial \alpha_1\partial \theta}=0=w_{31}, w_{23}=\dfrac{\partial^2 w}{\partial \alpha_2\partial \theta}=0=w_{32}, w_{33}=0.
\end{align*}
From \eqref{loglike}, we find the following quantities pertaining to likelihood function
\begin{align*}
	\mathbb{L}_{11}&=-\dfrac{n}{\alpha_1^2},\mathbb{L}_{12}=0=\mathbb{L}_{21},\mathbb{L}_{22}=-\dfrac{m}{\alpha_2^2}, \mathbb{L}_{13}=\dfrac{1}{\theta}=\mathbb{L}_{31}, \mathbb{L}_{23}=\dfrac{1}{\theta}=\mathbb{L}_{32},\mathbb{L}_{33}=-\dfrac{\alpha_1+\alpha_2}{\theta}\\
	\mathbb{L}_{111}&=-\dfrac{2n}{\alpha_1^3},\mathbb{L}_{133}=-\dfrac{1}{\theta^2}=\mathbb{L}_{331}=\mathbb{L}_{313},\mathbb{L}_{222}=-\dfrac{2m}{\alpha_2^3}, \mathbb{L}_{233}=-\dfrac{1}{\theta^2}=\mathbb{L}_{332}=\mathbb{L}_{323},\\
	 \mathbb{L}_{333}&=-\dfrac{2(\alpha_1+\alpha_2)}{\theta}.
\end{align*}
After taking logarithm of \eqref{jointprior}, we obtain
\begin{equation*}
\rho_j=\dfrac{\partial \rho}{\partial \alpha_j}=\dfrac{a_j-1}{\alpha_j}-b_j,~j=1,2~\textup{and}~\rho_3=\dfrac{\partial \rho}{\partial \theta}=\dfrac{a_3-1}{\theta}-b_3.
\end{equation*}
Then the Bayes estimator under SEL is
\begin{equation}\label{lind_sel}
	\delta_{s}^{L}=w+\left(\sum_{i=1}^{2}w_id_i+d_4+d_5\right)+\dfrac{1}{2}\left(A\sum_{i=1}^{2}w_i\sigma_{1i}+B\sum_{i=1}^{2}w_i\sigma_{2i}+C\sum_{i=1}^{2}w_i\sigma_{3i}\right)
\end{equation}
Further, for the Bayes estimator under LINEX function, all the quantities defined above will remain same except $w(\alpha_1,\alpha_2,\theta)$ and its respective derivatives. So, $w(\alpha_1,\alpha_2,\theta)$ under LINEX function is
\begin{align*}
w(\alpha_1,\alpha_2,\theta)=\textup{exp}\left(-c\sum_{p=s}^{k}\sum_{u=0}^{k-p}(-1)^u{k \choose p}{k-p \choose u}\left(\dfrac{\alpha_2}{\alpha_1(p+u)+\alpha_2}\right)\right)
\end{align*}
Again \eqref{lind_sel} is used to obtain Bayes estimator under LINEX function i.e., $	\delta_{l}^{L}$ with modified $w(\alpha_1,\alpha_2,\theta)$ and its respective derivatives. Once again all quantities are evaluated at the MLEs of the unknown parameters.
\subsubsection{Markov Chain Monte Carlo Method}\label{MCMC}
\begin{figure}[ht] 
	\begin{subfigure}[b]{0.5\linewidth}
		\centering
		\includegraphics[width=0.85\linewidth]{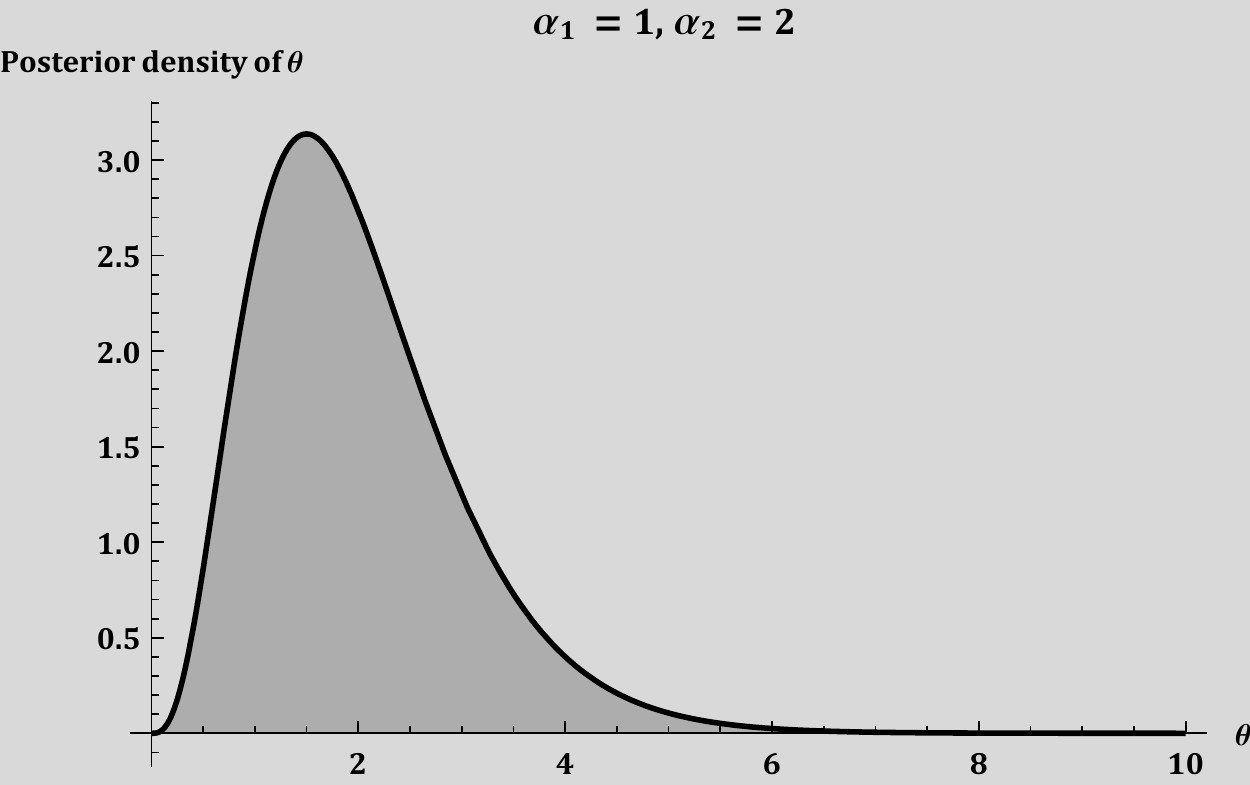} 
		\vspace{4ex}
	\end{subfigure}
	\begin{subfigure}[b]{0.5\linewidth}
		\centering
		\includegraphics[width=0.85\linewidth]{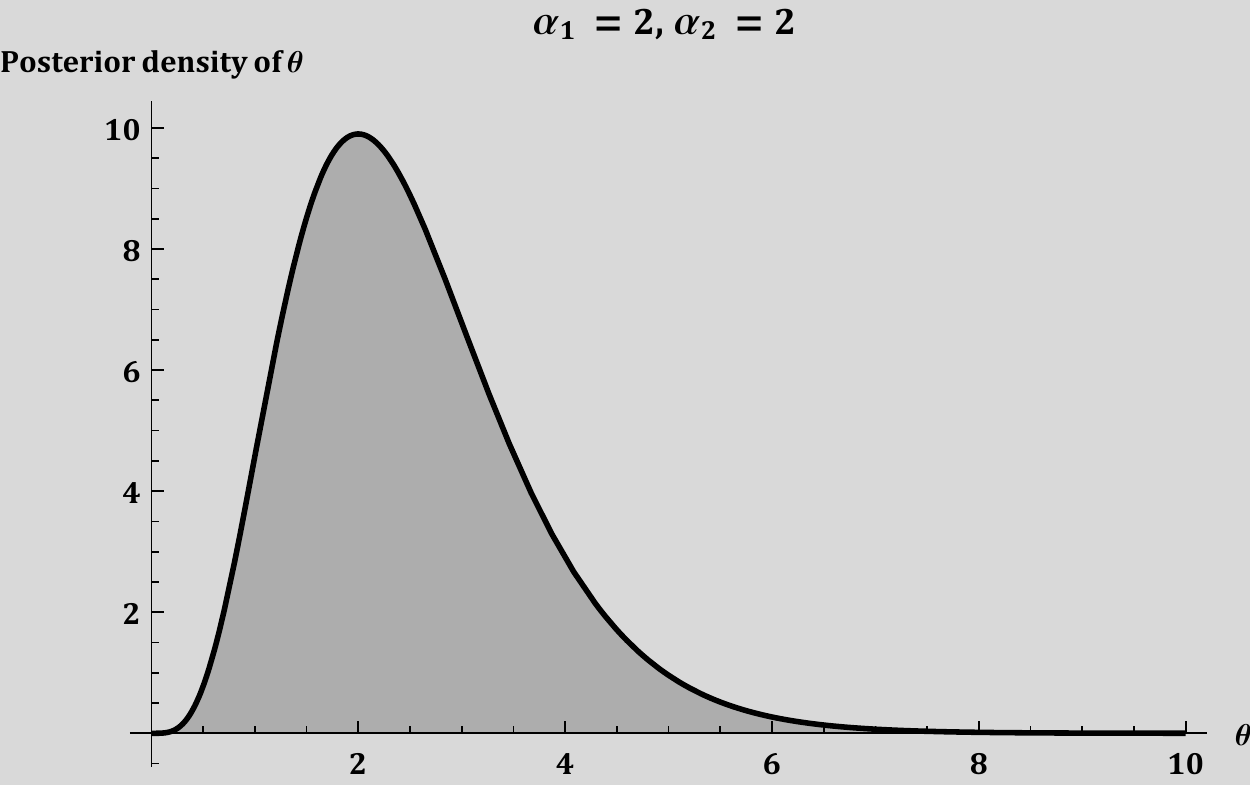} 
		\vspace{4ex}
	\end{subfigure} 
	\begin{subfigure}[b]{0.5\linewidth}
		\centering
		\includegraphics[width=0.85\linewidth]{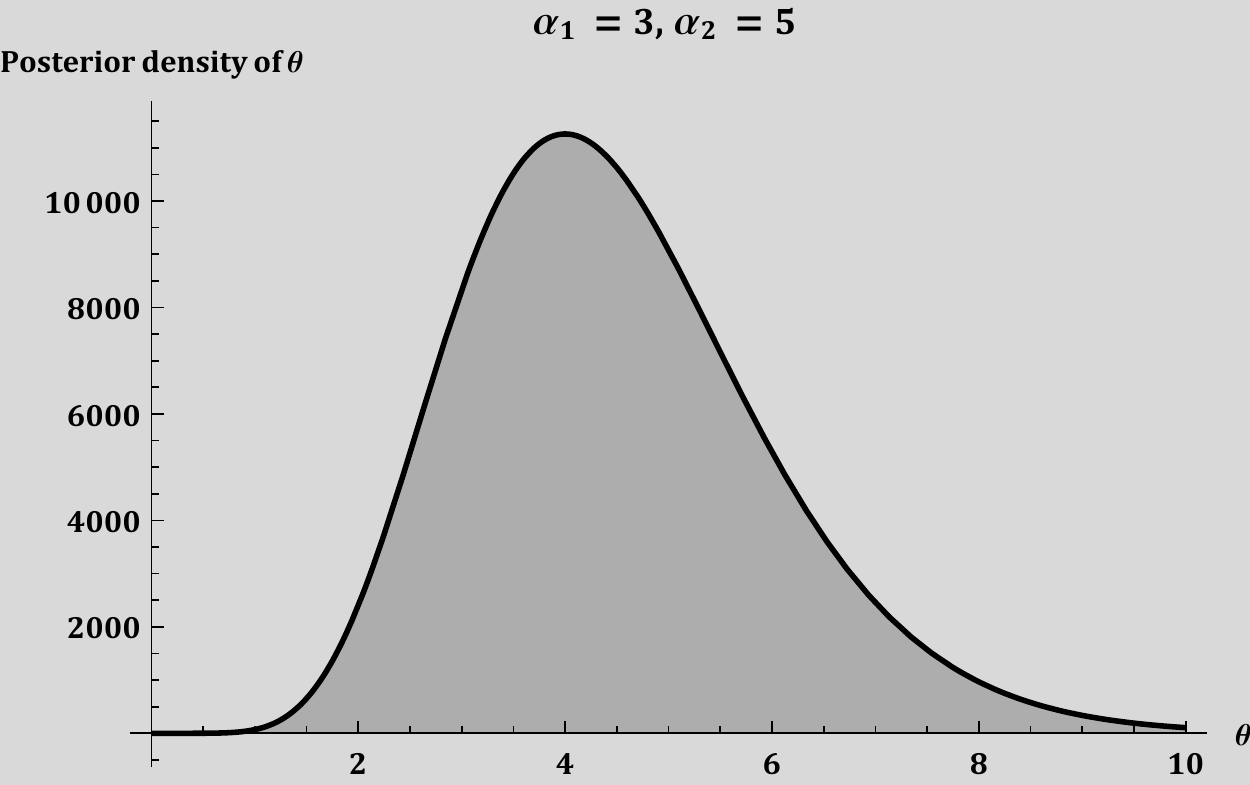} 
	\end{subfigure}
	\begin{subfigure}[b]{0.5\linewidth}
		\centering
		\includegraphics[width=0.85\linewidth]{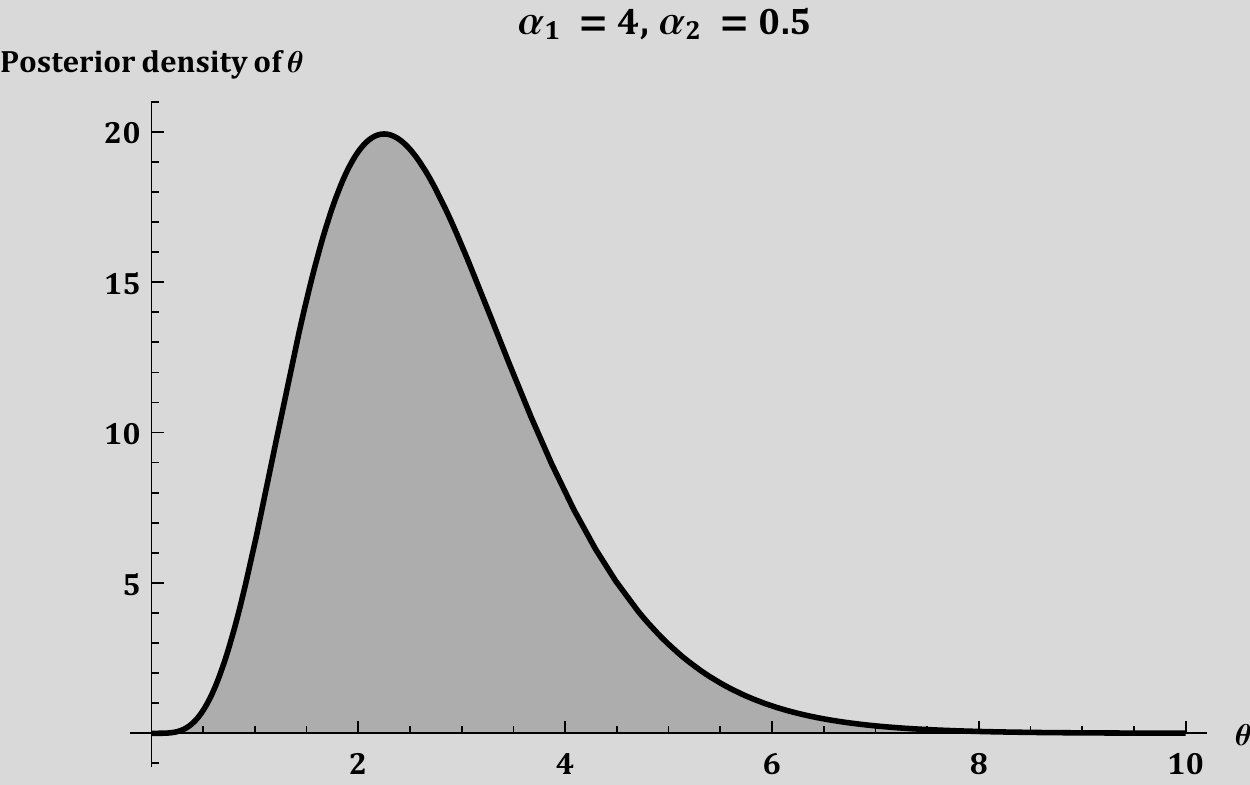} 
	\end{subfigure} 
	\caption{Various posterior densities of $\theta$}
	\label{plottheta} 
\end{figure}

The MCMC techniques is a powerful tool for obtaining approximate Bayes estimator with the help of marginal posterior densities.  
The joint posterior density of $\alpha_1,\alpha_2,\theta$ under the prior densities given in \eqref{jointprior}  is
\begin{equation*}
	\pi(\alpha_1,\alpha_2,\theta|data)\propto \alpha_1^{a_1-1}\alpha_2^{a_2-1}\theta^{\alpha_1+\alpha_2+a_3-1}e^{-\alpha_1(b_1-\ln r_n)}e^{-\alpha_2(b_2-\ln s_m)}e^{-b_3\theta}.
\end{equation*}
The marginal posterior densities of $\alpha_1$, $\alpha_2$ and $\theta$ are given as
\begin{align*}
	\pi(\alpha_1|\alpha_2,\theta,data)&\sim Gamma(n+a_1,b_1+\ln(r_n/\theta)),\qquad\alpha_1>0,\\
	\pi(\alpha_2|\alpha_1,\theta,data)&\sim Gamma(m+a_2,b_2+\ln(s_m/\theta)),\qquad\alpha_2>0,\\
	\pi(\theta|\alpha_1,\alpha_2,data)&\propto \theta^{\alpha_1+\alpha_2+a_3} e^{-b_3\theta},\qquad\qquad 0<\theta<\textup{min}\{R_1,S_1\}.
\end{align*}
We see that, marginal posterior densities of $\alpha_1$ and $\alpha_2$ have closed form of gamma distribution. Thus, using Gibbs sampling suggested by \cite{geman1987stochastic}, we generate random sample from marginal posterior densities of $\alpha_1$ and $\alpha_2$ whereas marginal posterior density of $\theta$ does not reduce to any analytical form of known distributions. But from Figure \ref{plottheta}, it can be seen that the marginal posterior density of $\theta$ is unimodal and has roughly symmetric form. So, we use Metropolis Hasting algorithm with proposal density of normal distribution suggested by \cite{gelman2013bayesian}. Hence, the following Metropolis-Hasting within Gibbs sampling algorithm is used
\begin{description}
	\item[(i)] Initiate with an initial value of $(\alpha_1^{(0)},\alpha_2^{(0)},\theta^{(0)}).$
	\item[(ii)] Set t=1.
	\item[(iii)] Generate $\theta^{(t)}$ by Metropolis-Hasting algorithm using normal distribution as proposal density. 
	\item[(iv)] Generate $\alpha_1^{(t)}$ from $Gamma\left (n+a_1,b_1+\ln\left(\tfrac{r_n}{\theta^{(t-1)}}\right)\right).$
	\item[(v)] Generate $\alpha_2^{(t)}$ from $Gamma\left (m+a_2,b_2+\ln\left(\tfrac{s_m}{\theta^{(t-1)}}\right)\right).$
	\item[(vi)] Compute $$R^{(t)}_{s;k}= \sum_{p=s}^{k}\sum_{u=0}^{k-p}(-1)^u{k \choose p}{k-p \choose u}\left(\dfrac{\alpha_2^{(t)}}{\alpha_1^{(t)}(p+u)+\alpha_2^{(t)}}\right).$$
	\item[(vii)] Repeat (ii)-(vi) $t=1,2,\ldots \ldots,T$ times.
\end{description} 
The generated sample is used to obtain the Bayes estimates of \rsk under SEL and LINEX functions. Also, it is used to obtain the HPD credible intervals for \rsk by using the method proposed by \cite{chen1999monte}. The Bayes estimators under SEL and LINEX function respectively, are
\begin{equation}
	\delta_{s}^{MC}=\dfrac{1}{T}\left(\sum_{i=1}^{T}R_{s;k}^{(i)}\right) \qquad \textup{and} \qquad 	\delta_{l}^{MC}=-\dfrac{1}{c} \ln \left(\dfrac{1}{T}\sum_{i=1}^{T}e^{-cR_{s;k}^{(i)}}\right).
\end{equation}
 \section{Estimation of \rsk for known $\theta$ }\label{asymptotic}
 In this section we assume that the common scale parameter $\theta$ is known and therefore all derived results of this section are based on this assumption.
 \subsection{MLE and Asymptotic theory of \rsk} 
 For known $\theta$, the MLEs of $\alpha_1$ and $\alpha_2$ can be rewritten as
 \begin{equation}\label{mlealphaknown}
 	\hat{\alpha_1}=\dfrac{n}{\ln R_n-\ln \theta} \quad  \& \quad	\hat{\alpha_2}=\dfrac{m}{\ln S_m-\ln \theta}
 \end{equation}
 Thus, MLE of \rsk is the same as defined as in \eqref{mlersk} with the values of $\hat{\alpha}_1$ and $\hat{\alpha}_2$ given in \eqref{mlealphaknown}. Further, we observe that the closed form distribution of \rsk is quite difficult to obtain. Therefore, asymptotic behavior of \rsk is discussed. For this purpose, we consider the Fisher information matrix of $\boldsymbol{\alpha}  =(\alpha_1,\alpha_2)$ as 
  $$ A(\boldsymbol{\alpha}) = -\left[\begin{array}{cc}
E \left(\dfrac{\partial^2 \ln L}{\partial \alpha_1^2}\right) & E\left(\dfrac{\partial^2 \ln L}{\partial \alpha_1 \partial \alpha_2}\right) \\ 
 &\\
 E\left(\dfrac{\partial^2 \ln }{\partial \alpha_2 \partial \alpha_1 }\right) & E\left(\dfrac{\partial^2 \ln L}{\partial \alpha_2^2}\right)
 \end{array} \right]=\left[\begin{array}{cc}
\frac{n}{\alpha_1^2} & 0 \\ 
 0& \frac{m}{\alpha_2^2}
 \end{array} \right]$$
The asymptotic variance of MLEs of $\alpha_1$ and $\alpha_2$ can easily be obtain from $A(\boldsymbol{\alpha})$. Also, the asymptotic variance of estimate of \rsk which contains $\hat{\alpha_1}$ and $\hat{\alpha_2}$ is given by (See \cite{rao2015burr}) as 
$$AV(\hat{R}_{s;k})=\left(\dfrac{\partial R_{s;k}}{\partial \alpha_1}\right)^2\dfrac{\alpha_1^2}{n}+\left(\dfrac{\partial R_{s;k}}{\partial \alpha_2}\right)^2 \dfrac{\alpha_2^2}{m}.$$ 
Now, from the asymptotic property of MLE, we establish that
$$\dfrac{\hat{R}_{s;k}-R_{s;k}}{AV(\hat{R}_{s;k})} \overset{d}{\to} N(0,1) \qquad n \to \infty,~m \to \infty.$$  
where $N(0,1)$ stands for standard normal distribution. Now, we construct the $(1-\beta)100\%$ asymptotic confidence interval for $R_{s;k}$ as
  \begin{equation}
  	\left(\hat{R}_{s;k}-z_{1-\beta/2}AV(\hat{R}_{s;k}),\hat{R}_{s;k}+z_{1-\beta/2}AV(\hat{R}_{s;k})\right).
  \end{equation}
  \subsection{UMVUE of $R _{s;k}$}
    In this section, UMVUE of \rsk is derived when $\theta$ is known, WLOG, assume $\theta=1$. The term \rsk can be rewritten as
  \begin{equation}\label{rsk.varphi}
  R_{s;k}=\sum_{p=s}^{k}\sum_{u=0}^{k-p}(-1)^u{k \choose p}{k-p \choose u}\varphi(\alpha_1,\alpha_2),
  \end{equation}
  where
  \begin{equation}\label{varphi}
  \varphi(\alpha_1,\alpha_2)=	\left(\dfrac{\alpha_2}{\alpha_1(p+u)+\alpha_2}\right).
  \end{equation}
Clearly, \rsk is a linear function of $\varphi(\alpha_1,\alpha_2)$ and $\varphi(\alpha_1,\alpha_2)$ depends on the unknown parameters $\alpha_1$ and $\alpha_2$.  Thus, it is enough to find the UMVUE of $\varphi(\alpha_1,\alpha_2).$
  \begin{equation*}
  L(\alpha_1,\alpha_2|\br,\bs)=\dfrac{\alpha_1^n\alpha_2^m\theta^{\alpha_1+\alpha_2}}{r_n^{\alpha_1}s_m^{\alpha_2}\prod_{i=1}^{n}r_i \prod_{i=1}^{m}s_i},\qquad\alpha_1>0,\alpha_2>0.
  \end{equation*}
  Applying Neyman Factorization Theorem, we find that $( R_n, S_m)$ is a sufficient for $(\alpha_1,\alpha_2).$ Define $U_1=\ln R_n,$ $U_2=\ln S_m,$ $V_1=\ln R_1$ and $V_2=\ln S_1.$ It is easy to verify that $(U_1,U_2)$ is complete and sufficient statistics for $(\alpha_1,\alpha_2).$ Consider a function
  \begin{equation}\label{indicator}
  \phi(V_1,V_2)=\begin{cases}
  1,&~if~ V_1>(p+u)V_2\\0&~\textup{otherwise}
  \end{cases}.
  \end{equation}
  It can be seen that $\phi(V_1,V_2)$ is an unbiased estimator of $\varphi(\alpha_1,\alpha_2).$ Now, using Rao-Blackwell and Lehman-Scheffe's Theorems, we get the UMVUE of $\varphi(\alpha_1,\alpha_2)$ is
  \begin{align}\label{varphi.UM}
  \varphi_{UM}(u_1,u_2)&=E\left[\phi(V_1,V_2)|U_1=u_1,U_2=u_2\right]\nonumber\\
  &=\iint\limits_\chi\phi(v_1,v_2) f_{V_1|U_1}(v_1|u_1)f_{V_2|U_2}(v_2|u_2)dv_1dv_2,
  \end{align}
  where $f_{V_1|U_1}$ and $f_{V_2|U_2}$ are the conditional \textit{pdfs} of $V_1|U_1$ and $V_2|U_2$ respectively. Now, we derive these conditional \textit{pdfs}. Using the formula of the joint \textit{pdf} of two upper records given in the book of \cite{ahsanullah1995record},
  the joint density of $(V_1,U_1)$ is
  \begin{equation}\label{joint.u1.v1}
  f_{V_1,U_1}(v_1,u_1)=\dfrac{\alpha_1^n}{(n-2)!}\left(u_1-v_1\right)^{n-2}e^{-u_1\alpha_1},\qquad 0<v_1<u_1<\infty.
  \end{equation}
  From \eqref{joint.u1.v1},  the marginal density of $U_1$ is
  \begin{equation}\label{marginal.u1}
  f_{U_1}(u_1)=\dfrac{\alpha_1^n}{(n-1)!}u_1^{n-1}e^{-u_1\alpha_1},\qquad u_1>0.
  \end{equation}
  Now, for fixed $u_1\in (0,\infty),$ the conditional distribution of $V_1$ given $U_1$ is
  \begin{equation}\label{conditional.v1.u1}
  f_{V_1|U_1}(v_1|u_1)=\dfrac{(n-1)(u_1-v_1)^{n-2}}{u_1^{n-1}},\qquad 0<v_1<u_1.
  \end{equation}
  Similarly, for fixed $u_2\in (0,\infty),$ the conditional distribution of $V_2$ given $U_2$ is
  \begin{equation}\label{conditional.v2.u2}
  f_{V_2|U_2}(v_2|u_2)=\dfrac{(m-1)(u_2-v_2)^{m-2}}{u_2^{m-1}},\qquad 0<v_2<u_2.
  \end{equation}
  Further, using \eqref{indicator}, \eqref{conditional.v1.u1} and \eqref{conditional.v2.u2} in \eqref{varphi.UM}, we get
  \begin{equation}
  \varphi_{UM}(u_1,u_2)=\iint\limits_\chi \dfrac{(n-1)(u_1-v_1)^{n-2}}{u_1^{n-1}} \dfrac{(m-1)(u_2-v_2)^{m-2}}{u_2^{m-1}} dv_1dv_2,
  \end{equation}
  where 	$\chi=\left\{ (v_1,v_2):0<v_1<u_1,0<v_2<u_2,(p+u)v_2<v_1\right\}.$ Depending on values of $u_1$ and $u_2$, the following two cases arise:
  \begin{description}
  	\item[Case (i)] When $0<u_2<\dfrac{u_1}{(p+u)};$ the UMVUE  is given as 
  	\begin{align*}
  	\varphi_{UM}(u_1,u_2)&=\dfrac{(n-1)(m-1)}{u_1^{n-1}u_2^{m-1}}\int_{0}^{u_2}\int_{(p+u)v_2}^{u_1}(u_1-v_1)^{n-2}(u_2-v_2)^{m-2}dv_1dv_2\\
  	&=\sum_{z_1=0}^{n-1}(-1)^{z_1}\dfrac{\Gamma m\Gamma n}{\Gamma(m+z_1) \Gamma(n-z_1)}\left(\dfrac{(p+u)u_2}{u_1}\right)^{z_1}
  	\end{align*}
  	\item[Case (ii)] When $0<\dfrac{u_1}{(p+u)}<u_2;$ the UMVUE  is given as 
  	\begin{align*}
  	\varphi_{UM}(u_1,u_2)&=\dfrac{(n-1)(m-1)}{u_1^{n-1}u_2^{m-1}}\int_{0}^{u_1/(p+u)}\int_{(p+u)v_2}^{u_1}(u_1-v_1)^{n-2}(u_2-v_2)^{m-2}dv_1dv_2\\
  	&=\sum_{z_1=0}^{m-2} (-1)^{z_1} \dfrac{\Gamma(m)\Gamma n}{\Gamma(n+z_1+1)\Gamma(m-z_1-1)} \left(\frac{u_1}{(p+u)u_2}\right)^{z_1+1}
  	\end{align*}
  \end{description}
  Thus, the UMVUE of $\varphi(\alpha_1,\alpha_2)$ is
  \begin{equation}\label{umvue}
  \varphi_{UM}(u_1,u_2)=\begin{cases}\sum_{z_1=0}^{n-1}(-1)^{z_1}\dfrac{\Gamma m\Gamma n}{\Gamma(m+z_1) \Gamma(n-z_1)}\left(\dfrac{(p+u)u_2}{u_1}\right)^{z_1},&  \textup{if}\quad 0<u_2<\dfrac{u_1}{(p+u)}\\
  \sum_{z_1=0}^{m-2} (-1)^{z_1} \dfrac{\Gamma(m)\Gamma n}{\Gamma(n+z_1+1)\Gamma(m-z_1-1)} \left(\frac{u_1}{(p+u)u_2}\right)^{z_1+1},&\textup{if}\quad 0<\dfrac{u_1}{(p+u)}<u_2.\end{cases}
  \end{equation}
It follows from \eqref{rsk.varphi} and \eqref{umvue} that the UMVUE of \rsk  is given by
  \begin{equation}\label{rsk.umvue}
  R^{UM}_{s;k}=\sum_{p=s}^{k}\sum_{u=0}^{k-p}(-1)^u{k \choose p}{k-p \choose u}\varphi_{UM}(u_1,u_2).
  \end{equation}
  
  \subsection{Bootstrapping}\label{bootstrapping}
  In this subsection, a well known parametric Bootstrap method is used for obtaining the various confidence intervals  for \rsk. For this purpose, we first generate the Bootstrap samples of \rsk as suggested by \cite{efron1994introduction}. The following algorithm is considered for generating the Bootstrap samples of \rsk.
  \begin{description}
  	\item[(i)] For primary samples of upper records, compute ML estimates of parameters as $\hat{\alpha}_1$ and $\hat{\alpha}_2$.
  	\item[(ii)] Generate independent Bootstrap samples of upper records from the Pareto distribution for known $\theta$ using parameters values as $\hat{\alpha}_1$ and $\hat{\alpha}_2$. After that, compute Bootstrap estimates of $\alpha_1$ and $\alpha_2$ as $\alpha_1^\ast$ and $\alpha_2^\ast,$ respectively.
  	\item[(iii)] Using the Bootstrap estimates $\alpha_1^\ast$ and $\alpha_2^\ast$, calculate $R_{s;k}^\ast$.
  	\item[(iv)] Replicate step-(ii) and step-(iii) \textit{B} times to obtain a sample of Bootstrap estimates $R_{s;k}^{1\ast},R_{s;k}^{2\ast},$ $\ldots, R_{s;k}^{B\ast}$ of \rsk.
  \end{description}
  Using the above generated Bootstrap sample of estimates of $R_{s;k},$ we find three types of confidence intervals of \rsk as follows:
  \begin{description}
  	\item[(I) \textit{Standard normal interval}:] This method provides the simplest confidence interval of $100(1-\beta)\%$ for \rsk using standard normal approach. The confidence interval for this approach is
  	\begin{equation}\label{bootnormal}
  		\left(\hat{R}_{s;k}-z_{1-\beta/2}\hat{\vartheta}^\ast,~\hat{R}_{s;k}+z_{1-\beta/2}\hat{\vartheta}^\ast \right).
  	\end{equation}
  	where $\hat{\vartheta}^\ast$ is estimate of the standard error based on $R_{s;k}^{1\ast},R_{s;k}^{2\ast},$ $\ldots, R_{s;k}^{B\ast}$ and $\hat{R}_{s;k}$ is the MLE of \rsk.
  	\item[(II) \textit{Percentile Bootstrap (Boot-p) interval}:] This method uses distribution function as a tool for obtaining the percentiles of any sample. Let $G(x)=P(R_{s;k}^\ast\leq x)$ be the distribution function of Bootstrap samples of \rsk for a given $x.$ Then $100(1-\beta)\%$ confidence interval of \rsk is
  	\begin{equation}
  		\left(G^{-1}\left(\frac{\beta}{2}\right),~ G^{-1}\left(1-\frac{\beta}{2}\right)\right),
  	\end{equation}
  	where $G^{-1}(t)$ is the solution of $G(x)=t.$
\item[(III) \textit{Bootstrap (Boot-t) interval}:] Let 
\begin{equation*}
	T_b^\ast=\dfrac{R_{s;k}^{b\ast}-\hat{R}}{\sqrt{\hat{var}\left(R_{s;k}^{b\ast}\right)}},\qquad\quad b=1,2,\ldots,B,
\end{equation*}
where $\sqrt{\hat{var}\left(R_{s;k}^{b\ast}\right)}$ is an estimate of standard error of $R_{s;k}^{b\ast}.$ For a large sample size, $\sqrt{\hat{var}\left(R_{s;k}^{b\ast}\right)}$ can be replaced by the asymptotic standard error  $\hat{\vartheta}^\ast$, used in \eqref{bootnormal}. If $T(x)=P(T_b^\ast\leq x)$ denote the distribution function of $T_b^\ast,$ then $100(1-\beta)\%$ confidence interval of \rsk is given by
\begin{equation}
		\left(\hat{R}-t^\ast_{1-\frac{\beta}{2}}\hat{\vartheta}^\ast,~ \hat{R}- t^\ast_{\frac{\beta}{2}}\hat{\vartheta}^\ast\right),
\end{equation}
where $t^\ast_{q}$ denotes $q^{th}$ quantile of $T_1^{\ast},T_2^{\ast},\ldots,T_B^{\ast}.$
  \end{description}
  \subsection{Bayesian Inference of \rsk}
  In this subsection, we derive the approximate Bayes estimators for known $\theta.$
\subsubsection{Lindley Approximation}
For the two parameter $(\alpha_1,\alpha_2)$, equation \eqref{lindley} reduces to
\begin{equation}\label{lindleyknown}
	E(w(\alpha_1,\alpha_2))=w+\left(w_1\tau_1+w_2\tau_2+\tau_3\right)+\dfrac{1}{2}\left[Q_1(w_1\sigma_{11}+w_2\sigma_{12})+Q_2(w_1\sigma_{21}+w_2\sigma_{22})\right]
\end{equation}
  where
  \begin{align*}
  \tau_i=&\rho_1\sigma_{i1}+\rho_2\sigma_{i2},\quad i=1,2,\\
  \tau_3=&\dfrac{1}{2}\left(w_{11}\sigma_{11}+w_{12}\sigma_{12}+w_{21}\sigma_{21}+w_{22}\sigma_{22}\right),\\
  Q_1=&\mathbb{L}_{111}\sigma_{11}+\mathbb{L}_{121}\sigma_{12}+\mathbb{L}_{211}\sigma_{21}+\mathbb{L}_{221}\sigma_{22},\\
  Q_2=&\mathbb{L}_{112}\sigma_{11}+\mathbb{L}_{122}\sigma_{12}+\mathbb{L}_{212}\sigma_{21}+\mathbb{L}_{222}\sigma_{22}.
\end{align*}
All expressions for the calculation of Bayes estimator can be calculated in same manner as in Section  \ref{Lindleymethod}.
\subsubsection{Mrakov Chain Monte Carlo Method}
It is evident from the Section \ref{MCMC} that, for known $\theta$, the marginal posterior densities of $\alpha_1$ and $\alpha_2$ follow gamma distribution. That is,
\begin{align*}
		\pi(\alpha_1|\alpha_2,data)&\sim Gamma(n+a_1,b_1+\ln(r_n/\theta))\\
	\pi(\alpha_2|\alpha_1,data)&\sim Gamma(m+a_2,b_2+\ln(s_m/\theta)).
\end{align*}
Now, using Gibbs algorithm, we generate data from these marginal densities. The algorithm is as follows
\begin{description}
	\item[(i)] Set t=1.
	\item[(ii)] Generate $\alpha_1^{(t)}$ from $\Gamma\left (n+a_1,b_1+\ln\left(\tfrac{r_n}{\theta}\right)\right).$
	\item[(iii)] Generate $\alpha_2^{(t)}$ from $\Gamma\left (m+a_2,b_2+\ln\left(\tfrac{s_m}{\theta}\right)\right).$
	\item[(iv)] Compute $$R^{(t)}_{s;k}= \sum_{p=s}^{k}\sum_{u=0}^{k-p}(-1)^u{k \choose p}{k-p \choose u}\left(\dfrac{\alpha_2^{(t)}}{\alpha_1^{(t)}(p+u)+\alpha_2^{(t)}}\right).$$
	\item[(v)] Repeat (i)-(iv) $t=1,2,\ldots\ldots,T$ times.
\end{description} 
The generated sample is used to obtain the Bayes estimates of \rsk under SEL and LINEX functions. Also, it is used to obtain the HPD credible intervals for \rsk by using the method proposed by \cite{chen1999monte}. The Bayes estimators under SEL and LINEX loss functions, respectively, are
\begin{equation}
\delta_{s}^{MC}=\dfrac{1}{T}\left(\sum_{i=1}^{T}R_{s;k}^{(i)}\right) \qquad \textup{and} \qquad 	\delta_{l}^{MC}=-\dfrac{1}{c} \ln \left(\dfrac{1}{T}\sum_{i=1}^{T}e^{-cR_{s;k}^{(i)}}\right).
\end{equation}
\section{Simulation Study}
In this section, a simulation study is conducted to exhibit the performances of the various estimator derived in this article. For this purpose, various sample sizes of upper record values are generated from Pareto distribution using the method discussed by \cite{wang2013reliability}. For performances of the estimator a Monte Carlo study is conducted for 1000 replications and results are reported in the following tables in terms of average estimate (AE), Mean squared error (MSE), coverage probability (CP) and average confidence interval length (AL). The results are presented in two cases: Case (I) discusses the results for unknown $\theta$ and Case (II) discusses the results for known $\theta.$ The following algorithm is used for calculation of results reported in the Tables.
\begin{enumerate}
	\item Generate samples of upper record values from $P(\alpha_1,\theta)$ and $P(\alpha_2,\theta)$ using some predefined values of parameters.
	\item Using the generated samples and applying the method discussed in previous sections, calculate estimator $\tilde{R}_{s;k}.$
	\item Repeat above steps 1000 times and obtain ${\tilde{R}_{s;k}}^1,{\tilde{R}_{s;k}}^2,\ldots,{\tilde{R}_{s;k}}^{1000}.$
	\item Now use the following equations to calculate the AE and MSE
	\begin{equation*}
		\textup{AE}=\dfrac{1}{1000}\sum_{i=1}^{1000}{\tilde{R}_{s;k}}^i\qquad\&\qquad 	\textup{MSE}=\dfrac{1}{1000}\sum_{i=1}^{1000}\left({\tilde{R}{s;k}}^i-R_{s;k}\right)^2,
	\end{equation*}
	where $R_{s;k}$ is true value.
\end{enumerate}
Also,the coverage probability can be obtain by using the following algorithm
\begin{enumerate}
	\item Generate samples of upper record values from $P(\alpha_1,\theta)$ and $P(\alpha_2,\theta)$ using some predefined values of parameters.
	\item Calculate the confidence interval for $R_{s;k}.$
	\item Repeat above steps, $1000$ times and find the number of intervals $(p)$ containing $R_{s;k}.$
	\item Now the coverage probability is $\tfrac{p}{1000}$.
\end{enumerate}
\section*{Case I when $\theta$ is unknown}
In this case, we consider 2-out-of-(4,5,6) and 3-out-of-(4,5,6) components i.e., $(s,k)=\{(2,4),(2,5),$ $(2,6),(3,4),(3,5),(3,6)\}$  for the numerical aspects. The values of \rsk for given set of $(s,k)$ are $R_{2;4}=0.80
$, $R_{2;5}=0.8571$, $R_{2;6}=0.8929$, $R_{3;4}=0.60$, $R_{3;5}=0.7143$ and $R_{3;6}=0.7857$ for $(\alpha_1,\alpha_2,\theta)=(2,4,1.5).$ All results reported in Table [\ref{unknown.A.and.mse.theta1.5.prior1}-\ref{unknown.hpd.alpha2.4}] are for unknown value of $\theta$.
 AEs and MSEs of the derived estimators are reported in Table[\ref{unknown.A.and.mse.theta1.5.prior1}-\ref{unknown.A.and.mse.theta1.5.prior2.continued}] with different sets of unknown quantities.
 Table[\ref{unknown.A.and.mse.theta1.5.prior1}-\ref{unknown.A.and.mse.theta1.5.prior1.continued}] reports AEs and MSEs of estimators of \rsk for prior $(a_1,a_2,a_3)=(2,2,2)$ and $(b_1,b_2,b_3)=(1.5,1.5,1.5)$. From these Tables, one can observe that Bayes estimator of Lindley method provides smaller MSEs in comparison with MCMC Bayes estimators for SEL and LINEX loss functions. 
The author also observe from  Table[\ref{unknown.A.and.mse.theta1.5.prior2}-\ref{unknown.A.and.mse.theta1.5.prior2.continued}] that for prior  $(a_1,a_2,a_3)=(3,3,3)$ and $(b_1,b_2,b_3)=(1.5,1.5,1.5),$ MCMC Bayes estimator shows smaller MSE for SEL function whereas for LINEX loss function Lindley Bayes estimator shows smaller MSE.. 
 
\begin{figure}[htbp]
	\centering
	\includegraphics[width=0.9\linewidth]{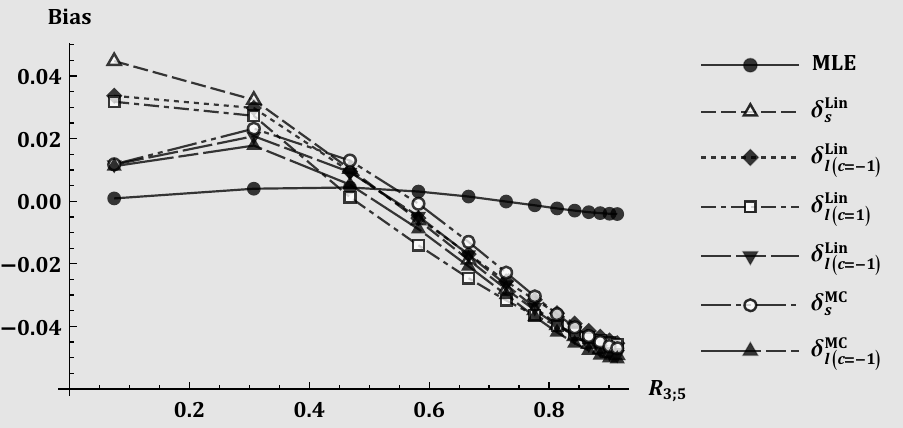}
	\caption{Bias plot of estimators when $\theta$ is unknown for $R_{3;5}$}
	\label{fig:bias1}
\end{figure}

\begin{figure}[htbp]
	\centering
	\includegraphics[width=0.9\linewidth]{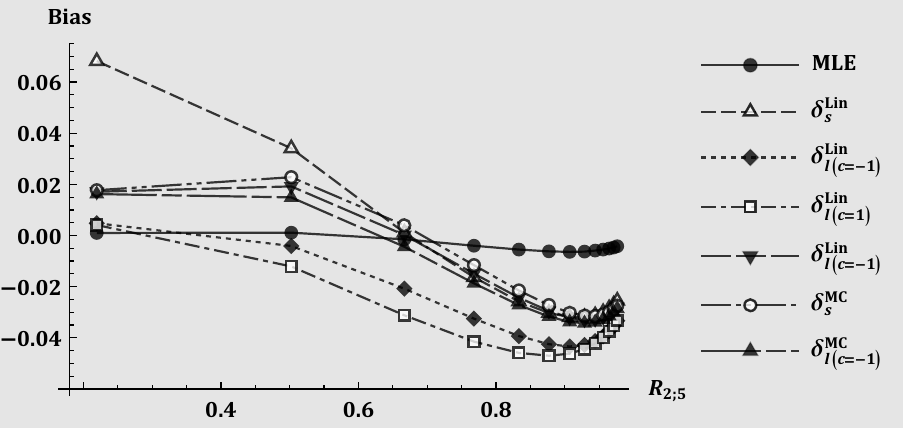}
	\caption{Bias plot of estimators when $\theta$ is unknown for $R_{2;5}$}
	\label{fig:bias3}
\end{figure}

    The HPD intervals are also obtained by using the method of \cite{chen1999monte} and reported in Table [\ref{unknown.hpd.alpha2.3}-\ref{unknown.hpd.alpha2.4} ] with CPs, ALs. From these tables we observe that, as we increase the sample sizes $n$ and $m$, the CPs are tending towards the desired level of significance. Also, the ALs are decreasing for large $n$ and $m$. The average Biases of the estimators are exhibited in Figure [\ref{fig:bias1}-\ref{fig:bias3}]. The biases are calculated with respective to the MSS reliability ranging from 0.1 to 0.9 as depicted in the Plots. Figure [\ref{fig:bias1}] shows the behavior of bias obtained at $(s,k)=(3,5)$, $n=20$, $m=20$ and $\theta=2$ for MSS reliability ranging from 0.1 to 0.9. Figure [\ref{fig:bias3}] shows the behavior of bias calculated at $(s,k)=(2,5)$,  $n=20$, $m=20$ and $\theta=2$. It is observed from Figure [\ref{fig:bias1}], that  bias of Bayes estimators has negative as well as positive values whereas bias of MLE remains close to 0 for increasing reliability. Also, Figure [\ref{fig:bias3}] shows significant variation among Bayes estimators whereas bias of MLE again remains close to 0  for increasing reliability. 
\begin{table}[htbp]
	\centering
	\caption{For $(\alpha_1,\alpha_2,\theta)=(2,4,1.5)$ and Prior: $(a_1,a_2,a_3)=(2,2,2)$,$(b_1,b_2,b_3)=(1.5,1.5,1.5).$ }
	\begin{tabular}{crrrrrrrrr}
		\toprule
		\multirow{4}[8]{*}{$(n,m)$} & \multicolumn{1}{c}{\multirow{4}[8]{*}{$(s,k)$}} & \multicolumn{2}{c}{$\hat{R}_{s;k}$} & \multicolumn{6}{c}{ Bayes Lindley Method} \\
		\cmidrule{3-10}          &       & \multicolumn{1}{c}{\multirow{3}[6]{*}{AE}} & \multicolumn{1}{c}{\multirow{3}[6]{*}{MSE}} & \multicolumn{2}{c}{SEL} & \multicolumn{4}{c}{LINEX} \\
		\cmidrule{5-10}          &       &       &       & \multicolumn{1}{c}{\multirow{2}[4]{*}{AE}} & \multicolumn{1}{c}{\multirow{2}[4]{*}{MSE}} & \multicolumn{2}{c}{c=-1} & \multicolumn{2}{c}{c=1} \\
		\cmidrule{7-10}          &       &       &       &       &       & \multicolumn{1}{c}{AE} & \multicolumn{1}{c}{MSE} & \multicolumn{1}{c}{AE} & \multicolumn{1}{c}{MSE} \\
		\midrule
		(10,10) & \multirow{6}[2]{*}{(2,4)} & \multicolumn{1}{c}{0.7949} & \multicolumn{1}{c}{0.0129} & 0.6955 & 0.0109 & \multicolumn{1}{c}{0.6287} & \multicolumn{1}{c}{0.0293} & 0.6264 & 0.0301 \\
		(10,15) &       & \multicolumn{1}{c}{0.7883} & \multicolumn{1}{c}{0.0112} & 0.7296 & 0.0050 & \multicolumn{1}{c}{0.7787} & \multicolumn{1}{c}{0.0150} & 0.7689 & 0.0210 \\
		(10,20) &       & \multicolumn{1}{c}{0.7852} & \multicolumn{1}{c}{0.0107} & 0.7435 & 0.0032 & \multicolumn{1}{c}{0.7741} & \multicolumn{1}{c}{0.0127} & 0.7635 & 0.0130 \\
		(15,15) &       & \multicolumn{1}{c}{0.7967} & \multicolumn{1}{c}{0.0087} & 0.7397 & 0.0036 & \multicolumn{1}{c}{0.7842} & \multicolumn{1}{c}{0.0092} & 0.7794 & 0.0124 \\
		(15,20) &       & \multicolumn{1}{c}{0.7929} & \multicolumn{1}{c}{0.0078} & 0.7441 & 0.0031 & \multicolumn{1}{c}{0.7648} & \multicolumn{1}{c}{0.0072} & 0.7614 & 0.0115 \\
		(20,20) &       & \multicolumn{1}{c}{0.7970} & \multicolumn{1}{c}{0.0065} & 0.6812 & 0.0038 & \multicolumn{1}{c}{0.7216} & \multicolumn{1}{c}{0.0061} & 0.7180 & 0.0103 \\ \midrule
		(10,10) &\multirow{6}[2]{*}{(2,5)} & \multicolumn{1}{c}{0.8469} & \multicolumn{1}{c}{0.0107} & 0.7638 & 0.0087 & \multicolumn{1}{c}{0.6931} & \multicolumn{1}{c}{0.0269} & 0.6924 & 0.0271 \\
		(10,15) &       & \multicolumn{1}{c}{0.8420} & \multicolumn{1}{c}{0.0094} & 0.7889 & 0.0047 & \multicolumn{1}{c}{0.8347} & \multicolumn{1}{c}{0.0250} & 0.8259 & 0.0210 \\
		(10,20) &       & \multicolumn{1}{c}{0.8396} & \multicolumn{1}{c}{0.0089} & 0.8018 & 0.0031 & \multicolumn{1}{c}{0.8294} & \multicolumn{1}{c}{0.0180} & 0.8206 & 0.0130 \\
		(15,15) &       & \multicolumn{1}{c}{0.8503} & \multicolumn{1}{c}{0.0071} & 0.8004 & 0.0030 & \multicolumn{1}{c}{0.8403} & \multicolumn{1}{c}{0.0113} & 0.8369 & 0.0114 \\
		(15,20) &       & \multicolumn{1}{c}{0.8474} & \multicolumn{1}{c}{0.0064} & 0.8037 & 0.0029 & \multicolumn{1}{c}{0.8233} & \multicolumn{1}{c}{0.0111} & 0.8211 & 0.0113 \\
		(20,20) &       & \multicolumn{1}{c}{0.8516} & \multicolumn{1}{c}{0.0052} & 0.7445 & 0.0027 & \multicolumn{1}{c}{0.7836} & \multicolumn{1}{c}{0.0080} & 0.7810 & 0.0098 \\
		\midrule
		(10,10) & \multirow{6}[2]{*}{(2,6)} & \multicolumn{1}{c}{0.8796} & \multicolumn{1}{c}{0.0088} & 0.8104 & 0.0068 & \multicolumn{1}{c}{0.7387} & \multicolumn{1}{c}{0.0238} & 0.7391 & 0.0236 \\
		(10,15) &       & \multicolumn{1}{c}{0.8759} & \multicolumn{1}{c}{0.0078} & 0.8281 & 0.0042 & \multicolumn{1}{c}{0.8700} & \multicolumn{1}{c}{0.0150} & 0.8625 & 0.0190 \\
		(10,20) &       & \multicolumn{1}{c}{0.8740} & \multicolumn{1}{c}{0.0075} & 0.8399 & 0.0028 & \multicolumn{1}{c}{0.8645} & \multicolumn{1}{c}{0.0138} & 0.8574 & 0.0130 \\
		(15,15) &       & \multicolumn{1}{c}{0.8839} & \multicolumn{1}{c}{0.0057} & 0.8405 & 0.0027 & \multicolumn{1}{c}{0.8759} & \multicolumn{1}{c}{0.0098} & 0.8735 & 0.0124 \\
		(15,20) &       & \multicolumn{1}{c}{0.8817} & \multicolumn{1}{c}{0.0051} & 0.8429 & 0.0025 & \multicolumn{1}{c}{0.8612} & \multicolumn{1}{c}{0.0076} & 0.8598 & 0.0071 \\
		(20,20) &       & \multicolumn{1}{c}{0.8857} & \multicolumn{1}{c}{0.0042} & 0.7876 & 0.0011 & \multicolumn{1}{c}{0.8250} & \multicolumn{1}{c}{0.0046} & 0.8232 & 0.0049 \\
		\midrule
		(10,10) &\multirow{6}[2]{*}{(3,4)} & \multicolumn{1}{c}{0.6110} & \multicolumn{1}{c}{0.0179} & 0.4808 & 0.0142 & \multicolumn{1}{c}{0.4317} & 0.0283 & \multicolumn{1}{c}{0.4256} & \multicolumn{1}{c}{0.0304} \\
		(10,15) &       & \multicolumn{1}{c}{0.6000} & \multicolumn{1}{c}{0.0146} & 0.5335 & 0.0114 & \multicolumn{1}{c}{0.5835} & 0.0230 & \multicolumn{1}{c}{0.5736} & \multicolumn{1}{c}{0.0270} \\
		(10,20) &       & \multicolumn{1}{c}{0.5953} & \multicolumn{1}{c}{0.0136} & 0.5477 & 0.0087 & \multicolumn{1}{c}{0.5820} & 0.0213 & \multicolumn{1}{c}{0.5691} & \multicolumn{1}{c}{0.0110} \\
		(15,15) &       & \multicolumn{1}{c}{0.6076} & \multicolumn{1}{c}{0.0121} & 0.5388 & 0.0047 & \multicolumn{1}{c}{0.5901} & 0.0150 & \multicolumn{1}{c}{0.5821} & \multicolumn{1}{c}{0.0093} \\
		(15,20) &       & \multicolumn{1}{c}{0.6016} & \multicolumn{1}{c}{0.0106} & 0.5451 & 0.0044 & \multicolumn{1}{c}{0.5668} & 0.0111 & \multicolumn{1}{c}{0.5603} & \multicolumn{1}{c}{0.0086} \\
		(20,20) &       & \multicolumn{1}{c}{0.6050} & \multicolumn{1}{c}{0.0091} & 0.4800 & 0.0040 & \multicolumn{1}{c}{0.5193} & 0.0065 & \multicolumn{1}{c}{0.5137} & \multicolumn{1}{c}{0.0074} \\
		\midrule
		(10,10) & \multirow{6}[2]{*}{(3,5)} & 0.7169 & 0.0171 & 0.5931 & 0.0147 & 0.5323 & 0.0331 & 0.5272 & 0.0350 \\
		(10,15) &       & 0.7077 & 0.0145 & 0.6407 & 0.0084 & 0.6947 & 0.0240 & 0.6833 & 0.0170 \\
		(10,20) &       & 0.7036 & 0.0137 & 0.6560 & 0.0064 & 0.6913 & 0.0190 & 0.6776 & 0.0130 \\
		(15,15) &       & 0.7163 & 0.0117 & 0.6486 & 0.0043 & 0.7004 & 0.0120 & 0.6931 & 0.0120 \\
		(15,20) &       & 0.7111 & 0.0104 & 0.6546 & 0.0036 & 0.6772 & 0.0114 & 0.6717 & 0.0101 \\
		(20,20) &       & 0.7152 & 0.0088 & 0.5862 & 0.0024 & 0.6288 & 0.0073 & 0.6235 & 0.0083 \\
		\midrule
		(10,10) & \multirow{6}[2]{*}{(3,6)} & 0.7817 & 0.0151 & 0.6706 & 0.0133 & 0.6022 & 0.0337 & 0.5989 & 0.0349 \\
		(10,15) &       & 0.7743 & 0.0131 & 0.7104 & 0.0087 & 0.7641 & 0.0220 & 0.7527 & 0.0181 \\
		(10,20) &       & 0.7708 & 0.0125 & 0.7256 & 0.0056 & 0.7594 & 0.0170 & 0.7467 & 0.0120 \\
		(15,15) &       & 0.7832 & 0.0103 & 0.7202 & 0.0053 & 0.7693 & 0.0130 & 0.7634 & 0.0106 \\
		(15,20) &       & 0.7789 & 0.0093 & 0.7255 & 0.0046 & 0.7477 & 0.0114 & 0.7435 & 0.0080 \\
		(20,20) &       & 0.7833 & 0.0077 & 0.6581 & 0.0039 & 0.7010 & 0.0072 & 0.6966 & 0.0079 \\
		\bottomrule
	\end{tabular}%
\label{unknown.A.and.mse.theta1.5.prior1}%
\end{table}%

\begin{table}[htbp]
	\centering
	\caption{\textit{continued}}
	\begin{tabular}{rrrrrrrrrrr}
		\toprule
		\multicolumn{2}{c}{Bayes Lindley} &       & \multicolumn{8}{c}{Bayes MCMC } \\
		\midrule
		\multicolumn{2}{c}{LINEX} &       & \multicolumn{2}{c}{SEL} & \multicolumn{6}{c}{LINEX} \\
		\midrule
		\multicolumn{2}{c}{c=1.5} &       & \multicolumn{1}{c}{\multirow{2}[4]{*}{AE}} & \multicolumn{1}{c}{\multirow{2}[4]{*}{MSE}} & \multicolumn{2}{c}{c=-1} & \multicolumn{2}{c}{c=1} & \multicolumn{2}{c}{c=1.5} \\
		\cmidrule{1-3}\cmidrule{6-11}    \multicolumn{1}{c}{AE} & \multicolumn{1}{c}{MSE} &       &       &       & \multicolumn{1}{c}{AE} & \multicolumn{1}{c}{MSE} & \multicolumn{1}{c}{AE} & \multicolumn{1}{c}{MSE} & \multicolumn{1}{c}{AE} & \multicolumn{1}{c}{MSE} \\
		\midrule
		\multicolumn{1}{c}{0.6263} & \multicolumn{1}{c}{0.0302} &       & 0.7336 & 0.0110 & 0.7394 & 0.0101 & 0.7277 & 0.0121 & 0.7247 & 0.0126 \\
		\multicolumn{1}{c}{0.7662} & \multicolumn{1}{c}{0.0160} &       & 0.7565 & 0.0078 & 0.7611 & 0.0073 & 0.7519 & 0.0084 & 0.7495 & 0.0088 \\
		\multicolumn{1}{c}{0.7608} & \multicolumn{1}{c}{0.0150} &       & 0.7695 & 0.0064 & 0.7735 & 0.0060 & 0.7655 & 0.0069 & 0.7634 & 0.0071 \\
		\multicolumn{1}{c}{0.7783} & \multicolumn{1}{c}{0.0105} &       & 0.7495 & 0.0082 & 0.7535 & 0.0077 & 0.7454 & 0.0088 & 0.7433 & 0.0091 \\
		\multicolumn{1}{c}{0.7607} & \multicolumn{1}{c}{0.0095} &       & 0.7615 & 0.0065 & 0.7649 & 0.0062 & 0.7581 & 0.0069 & 0.7563 & 0.0071 \\
		\multicolumn{1}{c}{0.7173} & \multicolumn{1}{c}{0.0068} &       & 0.7584 & 0.0065 & 0.7614 & 0.0062 & 0.7553 & 0.0069 & 0.7537 & 0.0071 \\
		\cmidrule{1-2}\cmidrule{4-11}    \multicolumn{1}{c}{0.6926} & \multicolumn{1}{c}{0.0271} &       & 0.7903 & 0.0102 & 0.7954 & 0.0093 & 0.7850 & 0.0112 & 0.7823 & 0.0117 \\
		\multicolumn{1}{c}{0.8236} & \multicolumn{1}{c}{0.0011} &       & 0.8122 & 0.0070 & 0.8161 & 0.0065 & 0.8081 & 0.0076 & 0.8060 & 0.0079 \\
		\multicolumn{1}{c}{0.8184} & \multicolumn{1}{c}{0.0015} &       & 0.8244 & 0.0056 & 0.8277 & 0.0052 & 0.8209 & 0.0060 & 0.8192 & 0.0063 \\
		\multicolumn{1}{c}{0.8361} & \multicolumn{1}{c}{0.0004} &       & 0.8064 & 0.0074 & 0.8099 & 0.0069 & 0.8029 & 0.0080 & 0.8010 & 0.0082 \\
		\multicolumn{1}{c}{0.8207} & \multicolumn{1}{c}{0.0013} &       & 0.8180 & 0.0058 & 0.8209 & 0.0054 & 0.8150 & 0.0061 & 0.8135 & 0.0063 \\
		\multicolumn{1}{c}{0.7805} & \multicolumn{1}{c}{0.0059} &       & 0.8155 & 0.0058 & 0.8181 & 0.0054 & 0.8129 & 0.0061 & 0.8115 & 0.0063 \\
		\cmidrule{1-2}\cmidrule{4-11}    \multicolumn{1}{c}{0.7395} & \multicolumn{1}{c}{0.0235} &       & 0.8276 & 0.0092 & 0.8321 & 0.0084 & 0.8228 & 0.0101 & 0.8204 & 0.0106 \\
		\multicolumn{1}{c}{0.8605} & \multicolumn{1}{c}{0.0011} &       & 0.8482 & 0.0062 & 0.8516 & 0.0057 & 0.8447 & 0.0067 & 0.8429 & 0.0070 \\
		\multicolumn{1}{c}{0.8557} & \multicolumn{1}{c}{0.0014} &       & 0.8596 & 0.0048 & 0.8624 & 0.0045 & 0.8566 & 0.0052 & 0.8551 & 0.0054 \\
		\multicolumn{1}{c}{0.8730} & \multicolumn{1}{c}{0.0004} &       & 0.8435 & 0.0065 & 0.8465 & 0.0061 & 0.8404 & 0.0070 & 0.8388 & 0.0073 \\
		\multicolumn{1}{c}{0.8595} & \multicolumn{1}{c}{0.0011} &       & 0.8544 & 0.0050 & 0.8568 & 0.0047 & 0.8518 & 0.0053 & 0.8505 & 0.0055 \\
		\multicolumn{1}{c}{0.8228} & \multicolumn{1}{c}{0.0049} &       & 0.8525 & 0.0050 & 0.8547 & 0.0047 & 0.8502 & 0.0053 & 0.8491 & 0.0055 \\
		\cmidrule{1-2}\cmidrule{4-11}    \multicolumn{1}{c}{0.4246} & 0.0308 &       & 0.5438 & 0.0109 & 0.5503 & 0.0102 & 0.5373 & 0.0117 & 0.5340 & 0.0121 \\
		\multicolumn{1}{c}{0.5708} & 0.0009 &       & 0.5669 & 0.0085 & 0.5725 & 0.0081 & 0.5614 & 0.0089 & 0.5586 & 0.0092 \\
		\multicolumn{1}{c}{0.5658} & 0.0012 &       & 0.5805 & 0.0076 & 0.5855 & 0.0074 & 0.5755 & 0.0078 & 0.5729 & 0.0080 \\
		\multicolumn{1}{c}{0.5803} & 0.0004 &       & 0.5570 & 0.0086 & 0.5617 & 0.0082 & 0.5523 & 0.0090 & 0.5500 & 0.0093 \\
		\multicolumn{1}{c}{0.5590} & 0.0017 &       & 0.5690 & 0.0073 & 0.5732 & 0.0070 & 0.5649 & 0.0076 & 0.5628 & 0.0077 \\
		\multicolumn{1}{c}{0.5126} & 0.0076 &       & 0.5642 & 0.0071 & 0.5679 & 0.0068 & 0.5605 & 0.0074 & 0.5587 & 0.0075 \\
		\cmidrule{1-2}\cmidrule{4-11}    0.5266 & 0.0352 &       & 0.6486 & 0.0124 & 0.6554 & 0.0114 & 0.6416 & 0.0135 & 0.6381 & 0.0141 \\
		0.6801 & 0.0012 &       & 0.6731 & 0.0092 & 0.6787 & 0.0086 & 0.6673 & 0.0099 & 0.6644 & 0.0102 \\
		0.6742 & 0.0016 &       & 0.6872 & 0.0079 & 0.6922 & 0.0075 & 0.6821 & 0.0083 & 0.6795 & 0.0086 \\
		0.6915 & 0.0005 &       & 0.6641 & 0.0096 & 0.6689 & 0.0090 & 0.6591 & 0.0102 & 0.6566 & 0.0105 \\
		0.6705 & 0.0019 &       & 0.6769 & 0.0078 & 0.6811 & 0.0074 & 0.6726 & 0.0083 & 0.6704 & 0.0085 \\
		0.6224 & 0.0084 &       & 0.6726 & 0.0077 & 0.6764 & 0.0074 & 0.6688 & 0.0081 & 0.6669 & 0.0083 \\
		\cmidrule{1-2}\cmidrule{4-11}    0.5986 & 0.0350 &       & 0.7158 & 0.0125 & 0.7224 & 0.0114 & 0.7091 & 0.0137 & 0.7057 & 0.0143 \\
		0.7496 & 0.0013 &       & 0.7402 & 0.0090 & 0.7454 & 0.0083 & 0.7348 & 0.0097 & 0.7321 & 0.0101 \\
		0.7435 & 0.0018 &       & 0.7540 & 0.0074 & 0.7586 & 0.0070 & 0.7493 & 0.0080 & 0.7470 & 0.0082 \\
		0.7620 & 0.0006 &       & 0.7324 & 0.0094 & 0.7370 & 0.0088 & 0.7277 & 0.0101 & 0.7253 & 0.0104 \\
		0.7426 & 0.0019 &       & 0.7452 & 0.0075 & 0.7491 & 0.0071 & 0.7412 & 0.0080 & 0.7391 & 0.0083 \\
		0.6957 & 0.0081 &       & 0.7416 & 0.0075 & 0.7452 & 0.0071 & 0.7380 & 0.0079 & 0.7362 & 0.0082 \\
		\cmidrule{1-2}\cmidrule{4-11}    \end{tabular}%
	\label{unknown.A.and.mse.theta1.5.prior1.continued}%
\end{table}%

\begin{table}[htbp]
	\centering
	\caption{For $(\alpha_1,\alpha_2,\theta)=(2,4,1.5)$ and Prior: $(a_1,a_2,a_3)=(3,3,3)$,$(b_1,b_2,b_3)=(1.5,1.5,1.5).$ }
	\begin{tabular}{crrrrrrrrr}
		\toprule
		\multirow{4}[8]{*}{$(n,m)$} & \multicolumn{1}{c}{\multirow{4}[8]{*}{$(s,k)$}} & \multicolumn{8}{c}{ Bayes Lindley Method} \\
		\cmidrule{3-10}          &       & \multicolumn{2}{c}{SEL} & \multicolumn{6}{c}{LINEX} \\
		\cmidrule{3-10}          &       & \multicolumn{1}{c}{\multirow{2}[4]{*}{AE}} & \multicolumn{1}{c}{\multirow{2}[4]{*}{MSE}} & \multicolumn{2}{c}{c=-1} & \multicolumn{2}{c}{c=1} & \multicolumn{2}{c}{c=1.5} \\
		\cmidrule{5-10}          &       &       &       & \multicolumn{1}{c}{AE} & \multicolumn{1}{c}{MSE} & \multicolumn{1}{c}{AE} & \multicolumn{1}{c}{MSE} & \multicolumn{1}{c}{AE} & \multicolumn{1}{c}{MSE} \\
\midrule  
    (10,10) & \multirow{6}[2]{*}{(2,4)} & 0.6862 & 0.0186 & 0.6518 & 0.0220 & 0.6450 & 0.0240 & 0.6439 & 0.0244 \\
	(10,15) &       & 0.7469 & 0.0068 & 0.7546 & 0.0021 & 0.7435 & 0.0032 & 0.7407 & 0.0035 \\
	(10,20) &       & 0.7655 & 0.0050 & 0.7582 & 0.0017 & 0.7477 & 0.0027 & 0.7452 & 0.0030 \\
	(15,15) &       & 0.7335 & 0.0090 & 0.7802 & 0.0070 & 0.7757 & 0.0061 & 0.7747 & 0.0063 \\
	(15,20) &       & 0.7538 & 0.0066 & 0.7667 & 0.0011 & 0.7631 & 0.0014 & 0.7623 & 0.0014 \\
	(20,20) &       & 0.7512 & 0.0066 & 0.7286 & 0.0051 & 0.7243 & 0.0057 & 0.7234 & 0.0059 \\
	\midrule
	(10,10) & \multirow{6}[2]{*}{(2,5)} & 0.7527 & 0.0173 & 0.7144 & 0.0204 & 0.7098 & 0.0217 & 0.7090 & 0.0219 \\
	(10,15) &       & 0.8057 & 0.0062 & 0.8121 & 0.0020 & 0.8024 & 0.0030 & 0.8000 & 0.0033 \\
	(10,20) &       & 0.8220 & 0.0045 & 0.8148 & 0.0018 & 0.8063 & 0.0026 & 0.8043 & 0.0028 \\
	(15,15) &       & 0.7945 & 0.0081 & 0.8368 & 0.0054 & 0.8336 & 0.0071 & 0.8329 & 0.0073 \\
	(15,20) &       & 0.8120 & 0.0059 & 0.8249 & 0.0010 & 0.8226 & 0.0012 & 0.8221 & 0.0012 \\
	(20,20) &       & 0.8102 & 0.0060 & 0.7899 & 0.0045 & 0.7867 & 0.0050 & 0.7860 & 0.0051 \\
	\midrule
	(10,10) & \multirow{6}[2]{*}{(2,6)} & 0.7998 & 0.0137 & 0.7581 & 0.0182 & 0.7550 & 0.0190 & 0.7546 & 0.0191 \\
	(10,15) &       & 0.8442 & 0.0055 & 0.8491 & 0.0019 & 0.8410 & 0.0027 & 0.8389 & 0.0029 \\
	(10,20) &       & 0.8585 & 0.0039 & 0.8514 & 0.0017 & 0.8446 & 0.0023 & 0.8430 & 0.0025 \\
	(15,15) &       & 0.8346 & 0.0071 & 0.8729 & 0.0040 & 0.8707 & 0.0052 & 0.8702 & 0.0055 \\
	(15,20) &       & 0.8498 & 0.0052 & 0.8626 & 0.0009 & 0.8611 & 0.0010 & 0.8608 & 0.0010 \\
	(20,20) &       & 0.8486 & 0.0052 & 0.8305 & 0.0039 & 0.8282 & 0.0042 & 0.8277 & 0.0042 \\
	\midrule
	(10,10) & \multirow{6}[2]{*}{(3,4)} & 0.4704 & 0.0236 & 0.4560 & 0.0207 & 0.4455 & 0.0239 & 0.4436 & 0.0245 \\
	(10,15) &       & 0.5489 & 0.0067 & 0.5585 & 0.0017 & 0.5467 & 0.0028 & 0.5436 & 0.0032 \\
	(10,20) &       & 0.5720 & 0.0055 & 0.5645 & 0.0013 & 0.5515 & 0.0024 & 0.5483 & 0.0027 \\
	(15,15) &       & 0.5312 & 0.0091 & 0.5851 & 0.0060 & 0.5775 & 0.0066 & 0.5758 & 0.0073 \\
	(15,20) &       & 0.5573 & 0.0068 & 0.5691 & 0.0010 & 0.5624 & 0.0014 & 0.5610 & 0.0015 \\
	(20,20) &       & 0.5527 & 0.0069 & 0.5272 & 0.0053 & 0.5208 & 0.0063 & 0.5195 & 0.0065 \\
	\midrule
	(10,10) & \multirow{6}[2]{*}{(3,5)} & 0.5846 & 0.0234 & 0.5580 & 0.0244 & 0.5478 & 0.0277 & 0.5461 & 0.0283 \\
	(10,15) &       & 0.6584 & 0.0077     & 0.6683 & 0.0021 & 0.6551 & 0.0035 & 0.6516 & 0.0039 \\
	(10,20) &       & 0.6805 & 0.0060   	& 0.6733 & 0.0017 & 0.6596 & 0.0030 & 0.6564 & 0.0034 \\
	(15,15) &       & 0.6419 & 0.0104 & 0.6955 & 0.0081 & 0.6886 & 0.0077 & 0.6871 & 0.0081 \\
	(15,20) &       & 0.6665 & 0.0077 & 0.6795 & 0.0012 & 0.6737 & 0.0016 & 0.6725 & 0.0017 \\
	(20,20) &       & 0.6626 & 0.0078 & 0.6368 & 0.0060 & 0.6306 & 0.0070 & 0.6293 & 0.0072 \\
	\midrule
	(10,10) & \multirow{6}[2]{*}{(3,6)} & 0.6629 & 0.0217 & 0.6275 & 0.0250 & 0.6190 & 0.0278 & 0.6176 & 0.0283 \\
	(10,15) &       & 0.7290 & 0.0078 & 0.7381 & 0.0023 & 0.7252 & 0.0037 & 0.7219 & 0.0041 \\
	(10,20) &       & 0.7492 & 0.0058 & 0.7420 & 0.0019 & 0.7295 & 0.0032 & 0.7265 & 0.0035 \\
	(15,15) &       & 0.7143 & 0.0103 & 0.7648 & 0.0067 & 0.7592 & 0.0072 & 0.7580 & 0.0086 \\
	(15,20) &       & 0.7364 & 0.0076 & 0.7498 & 0.0013 & 0.7453 & 0.0016 & 0.7444 & 0.0017 \\
	(20,20) &       & 0.7334 & 0.0077 & 0.7087 & 0.0059 & 0.7034 & 0.0068 & 0.7023 & 0.0070 \\
	\bottomrule
\end{tabular}%
	\label{unknown.A.and.mse.theta1.5.prior2}%
\end{table}%
\begin{table}[htbp]
	\centering
	\caption{\textit{continued}}
	\begin{tabular}{crrrrrrrrr}
		\toprule
		\multirow{4}[8]{*}{$(n,m)$} & \multicolumn{1}{c}{\multirow{4}[8]{*}{$(s,k)$}} & \multicolumn{8}{c}{ Bayes MCMC Method} \\
		\cmidrule{3-10}          &       & \multicolumn{2}{c}{SEL} & \multicolumn{6}{c}{LINEX} \\
		\cmidrule{3-10}          &       & \multicolumn{1}{c}{\multirow{2}[4]{*}{AE}} & \multicolumn{1}{c}{\multirow{2}[4]{*}{MSE}} & \multicolumn{2}{c}{c=-1} & \multicolumn{2}{c}{c=1} & \multicolumn{2}{c}{c=1.5} \\
		\cmidrule{5-10}          &       &       &       & \multicolumn{1}{c}{AE} & \multicolumn{1}{c}{MSE} & \multicolumn{1}{c}{AE} & \multicolumn{1}{c}{MSE} & \multicolumn{1}{c}{AE} & \multicolumn{1}{c}{MSE} \\
		\midrule
		(10,10) & \multirow{6}[2]{*}{(2,4)} & 0.7344 & 0.0110 & 0.7398 & 0.0101 & 0.7290 & 0.0119 & 0.7262 & 0.0124 \\
		(10,15) &       & 0.7509 & 0.0085 & 0.7553 & 0.0080 & 0.7465 & 0.0092 & 0.7442 & 0.0095 \\
		(10,20) &       & 0.7605 & 0.0073 & 0.7644 & 0.0069 & 0.7566 & 0.0078 & 0.7546 & 0.0081 \\
		(15,15) &       & 0.7500 & 0.0082 & 0.7537 & 0.0077 & 0.7461 & 0.0088 & 0.7441 & 0.0090 \\
		(15,20) &       & 0.7585 & 0.0069 & 0.7617 & 0.0065 & 0.7551 & 0.0073 & 0.7534 & 0.0075 \\
		(20,20) &       & 0.7587 & 0.0065 & 0.7617 & 0.0062 & 0.7558 & 0.0068 & 0.7543 & 0.0070 \\
		\midrule
		(10,10) &\multirow{6}[2]{*}{(2,5)} & 0.7914 & 0.0101 & 0.7961 & 0.0093 & 0.7864 & 0.0110 & 0.7839 & 0.0115 \\
		(10,15) &       & 0.8072 & 0.0077 & 0.8110 & 0.0071 & 0.8033 & 0.0083 & 0.8013 & 0.0086 \\
		(10,20) &       & 0.8164 & 0.0065 & 0.8197 & 0.0061 & 0.8130 & 0.0070 & 0.8113 & 0.0072 \\
		(15,15) &       & 0.8070 & 0.0074 & 0.8103 & 0.0069 & 0.8036 & 0.0079 & 0.8019 & 0.0082 \\
		(15,20) &       & 0.8153 & 0.0061 & 0.8181 & 0.0057 & 0.8124 & 0.0065 & 0.8109 & 0.0067 \\
		(20,20) &       & 0.8160 & 0.0057 & 0.8185 & 0.0054 & 0.8134 & 0.0061 & 0.8121 & 0.0062 \\
		\midrule
		(10,10) & \multirow{6}[2]{*}{(2,6)} & 0.8287 & 0.0091 & 0.8329 & 0.0083 & 0.8244 & 0.0099 & 0.8221 & 0.0103 \\
		(10,15) &       & 0.8438 & 0.0068 & 0.8471 & 0.0063 & 0.8404 & 0.0073 & 0.8387 & 0.0076 \\
		(10,20) &       & 0.8524 & 0.0056 & 0.8552 & 0.0053 & 0.8495 & 0.0060 & 0.8481 & 0.0063 \\
		(15,15) &       & 0.8441 & 0.0065 & 0.8469 & 0.0061 & 0.8412 & 0.0070 & 0.8397 & 0.0072 \\
		(15,20) &       & 0.8520 & 0.0053 & 0.8543 & 0.0050 & 0.8495 & 0.0056 & 0.8483 & 0.0058 \\
		(20,20) &       & 0.8529 & 0.0050 & 0.8550 & 0.0047 & 0.8507 & 0.0053 & 0.8496 & 0.0054 \\ \midrule
		(10,10) & \multirow{6}[2]{*}{(3,4)} & 0.5440 & 0.0110 & 0.5500 & 0.0103 & 0.5379 & 0.0117 & 0.5349 & 0.0121 \\
		(10,15) &       & 0.5602 & 0.0091 & 0.5654 & 0.0087 & 0.5550 & 0.0095 & 0.5524 & 0.0098 \\
		(10,20) &       & 0.5699 & 0.0082 & 0.5746 & 0.0079 & 0.5651 & 0.0086 & 0.5628 & 0.0087 \\
		(15,15) &       & 0.5572 & 0.0087 & 0.5616 & 0.0083 & 0.5527 & 0.0091 & 0.5505 & 0.0093 \\
		(15,20) &       & 0.5653 & 0.0076 & 0.5693 & 0.0073 & 0.5613 & 0.0079 & 0.5594 & 0.0080 \\
		(20,20) &       & 0.5644 & 0.0071 & 0.5679 & 0.0069 & 0.5608 & 0.0074 & 0.5591 & 0.0075 \\
		\midrule
		(10,10) & \multirow{6}[2]{*}{(3,5)} & 0.6491 & 0.0124 & 0.6555 & 0.0115 & 0.6426 & 0.0134 & 0.6394 & 0.0140 \\
		(10,15) &       & 0.6665 & 0.0100 & 0.6718 & 0.0094 & 0.6610 & 0.0107 & 0.6583 & 0.0110 \\
		(10,20) &       & 0.6768 & 0.0088 & 0.6816 & 0.0083 & 0.6719 & 0.0093 & 0.6694 & 0.0096 \\
		(15,15) &       & 0.6644 & 0.0096 & 0.6690 & 0.0091 & 0.6597 & 0.0102 & 0.6573 & 0.0105 \\
		(15,20) &       & 0.6733 & 0.0082 & 0.6773 & 0.0078 & 0.6691 & 0.0087 & 0.6670 & 0.0089 \\
		(20,20) &       & 0.6729 & 0.0078 & 0.6766 & 0.0074 & 0.6692 & 0.0081 & 0.6674 & 0.0083 \\
		\midrule
		(10,10) &\multirow{6}[2]{*}{(3,6)} & 0.7166 & 0.0125 & 0.7227 & 0.0115 & 0.7104 & 0.0136 & 0.7072 & 0.0142 \\
		(10,15) &       & 0.7341 & 0.0098 & 0.7391 & 0.0091 & 0.7289 & 0.0105 & 0.7263 & 0.0109 \\
		(10,20) &       & 0.7443 & 0.0084 & 0.7488 & 0.0079 & 0.7398 & 0.0090 & 0.7375 & 0.0093 \\
		(15,15) &       & 0.7328 & 0.0094 & 0.7372 & 0.0088 & 0.7284 & 0.0101 & 0.7261 & 0.0104 \\
		(15,20) &       & 0.7418 & 0.0079 & 0.7456 & 0.0075 & 0.7380 & 0.0084 & 0.7360 & 0.0087 \\
		(20,20) &       & 0.7420 & 0.0075 & 0.7454 & 0.0071 & 0.7386 & 0.0079 & 0.7368 & 0.0081 \\
		\bottomrule
	\end{tabular}%
	\label{unknown.A.and.mse.theta1.5.prior2.continued}%
\end{table}%

 \begin{table}[htbp]
 	\begin{threeparttable}
 		\centering
 		\caption{For $(\alpha_1,\alpha_2,\theta)=(2,3,1.5)$, $Prior_1$: $(a_1,a_2,a_3)=(2,2,2)$,$(b_1,b_2,b_3)=(1.5,1.5,1.5)$ and $Prior_2$: $(a_1,a_2,a_3)=(1.5,1.5,1.5)$,$(b_1,b_2,b_3)=(2,2,2).$}
 		\begin{tabular}{cccccccccc}
 			\toprule
 			\multirow{3}[6]{*}{$(n,m)$} & \multirow{3}[6]{*}{$(s,k)$} & \multicolumn{4}{c}{$Prior_1$} & \multicolumn{4}{c}{$Prior_2$} \\
 			\cmidrule{3-10}          &       & \multicolumn{2}{c}{90\%} & \multicolumn{2}{c}{95\%} & \multicolumn{2}{c}{90\%} & \multicolumn{2}{c}{95\%} \\
 			\cmidrule{3-10}          &       & CP    & AL    & CP    & AL    & CP    & AL    & CP    & AL \\
 			\midrule
 			(10,10) & \multirow{6}[2]{*}{(2,4)} & 0.9280 & 0.3591 & 0.9780 & 0.4218 & 0.9540 & 0.3719 & 0.9860 & 0.4361 \\
 			(10,15) &       & 0.9140 & 0.3277 & 0.9720 & 0.3857 & 0.9480 & 0.3338 & 0.9880 & 0.3916 \\
 			(10,20) &       & 0.9120 & 0.3101 & 0.9660 & 0.3649 & 0.9300 & 0.3119 & 0.9700 & 0.3665 \\
 			(15,15) &       & 0.9100 & 0.3052 & 0.9620 & 0.3600 & 0.9220 & 0.3133 & 0.9760 & 0.3691 \\
 			(15,20) &       & 0.9140 & 0.2863 & 0.9700 & 0.3379 & 0.9360 & 0.2906 & 0.9740 & 0.3428 \\
 			(20,20) &       & 0.8960 & 0.2708 & 0.9420 & 0.3200 & 0.9160 & 0.2764 & 0.9680 & 0.3264 \\
 			\midrule
 			(10,10) & \multirow{6}[2]{*}{(2,6)} & 0.9280 & 0.3179 & 0.9780 & 0.3755 & 0.9580 & 0.3318 & 0.9880 & 0.3915 \\
 			(10,15) &       & 0.9120 & 0.2849 & 0.9800 & 0.3366 & 0.9420 & 0.2881 & 0.9900 & 0.3398 \\
 			(10,20) &       & 0.9140 & 0.2669 & 0.9580 & 0.3150 & 0.9220 & 0.2640 & 0.9680 & 0.3121 \\
 			(15,15) &       & 0.9120 & 0.2694 & 0.9600 & 0.3188 & 0.9320 & 0.2784 & 0.9780 & 0.3293 \\
 			(15,20) &       & 0.9160 & 0.2508 & 0.9660 & 0.2969 & 0.9300 & 0.2540 & 0.9720 & 0.3009 \\
 			(20,20) &       & 0.8900 & 0.2388 & 0.9480 & 0.2833 & 0.9180 & 0.2451 & 0.9640 & 0.2906 \\
 			\midrule
 			(10,10) & \multirow{6}[2]{*}{(3,4)} & 0.9300 & 0.3627 & 0.9740 & 0.4267 & 0.9480 & 0.3706 & 0.9840 & 0.4359 \\
 			(10,15) &       & 0.9340 & 0.3383 & 0.9740 & 0.3992 & 0.9560 & 0.3465 & 0.9920 & 0.408 \\
 			(10,20) &       & 0.9240 & 0.3243 & 0.9720 & 0.3830 & 0.9400 & 0.3319 & 0.9720 & 0.3913 \\
 			(15,15) &       & 0.8940 & 0.3090 & 0.9520 & 0.3652 & 0.9240 & 0.3141 & 0.9720 & 0.371 \\
 			(15,20) &       & 0.9140 & 0.2927 & 0.9620 & 0.3465 & 0.9340 & 0.2974 & 0.9740 & 0.3515 \\
 			(20,20) &       & 0.8880 & 0.2744 & 0.9440 & 0.3246 & 0.9020 & 0.2774 & 0.9640 & 0.3285 \\
 			\midrule
 			(10,10) & \multirow{6}[2]{*}{(3,6)} & 0.9240 & 0.3809 & 0.9780 & 0.4461 & 0.9500 & 0.3940 & 0.9860 & 0.4607 \\
 			(10,15) &       & 0.9140 & 0.3494 & 0.9720 & 0.4102 & 0.9480 & 0.3564 & 0.9880 & 0.4171 \\
 			(10,20) &       & 0.9080 & 0.3315 & 0.9660 & 0.3893 & 0.9280 & 0.3345 & 0.9700 & 0.392 \\
 			(15,15) &       & 0.9060 & 0.3250 & 0.9560 & 0.3825 & 0.9200 & 0.3333 & 0.9720 & 0.3917 \\
 			(15,20) &       & 0.9120 & 0.3057 & 0.9680 & 0.3601 & 0.9340 & 0.3105 & 0.9740 & 0.3654 \\
 			(20,20) &       & 0.8960 & 0.2890 & 0.9420 & 0.3409 & 0.9120 & 0.2946 & 0.9660 & 0.3473 \\
 			\bottomrule
 		\end{tabular}%
 		\label{unknown.hpd.alpha2.3}%
 	\end{threeparttable}
 \end{table}%

 \begin{table}[htbp]
 	\begin{threeparttable}
 		\centering
 		\caption{For $(\alpha_1,\alpha_2,\theta)$=$(2,4,1.5)$, $Prior_1$: $(a_1,a_2,a_3)$=$(2,2,2)$,$(b_1,b_2,b_3)$=$(1.5,1.5,1.5)$ and $Prior_2$: $(a_1,a_2,a_3)$=$(1.5,1.5,1.5)$,$(b_1,b_2,b_3)$=$(2,2,2).$}
 		\begin{tabular}{cccccccccc}
 			\toprule
 			\multirow{3}[6]{*}{$(n,m)$} & \multirow{3}[6]{*}{$(s,k)$} & \multicolumn{4}{c}{$Prior_1$} & \multicolumn{4}{c}{$Prior_2$} \\
 			\cmidrule{3-10}          &       & \multicolumn{2}{c}{90\%} & \multicolumn{2}{c}{95\%} & \multicolumn{2}{c}{90\%} & \multicolumn{2}{c}{95\%} \\
 			\cmidrule{3-10}          &       & CP    & AL    & CP    & AL    & CP    & AL    & CP    & AL \\
 			\midrule
 			(10,10) & \multirow{6}[2]{*}{(2,4)} & 0.9340 & 0.3359 & 0.9740 & 0.3955 & 0.9480 & 0.3514 & 0.9700 & 0.4134 \\
 			(10,15) &       & 0.9440 & 0.2993 & 0.9800 & 0.3528 & 0.9560 & 0.3059 & 0.9860 & 0.3606 \\
 			(10,20) &       & 0.9260 & 0.2787 & 0.9660 & 0.3288 & 0.9460 & 0.2801 & 0.9860 & 0.3305 \\
 			(15,15) &       & 0.9040 & 0.2811 & 0.9660 & 0.3322 & 0.9180 & 0.2911 & 0.9700 & 0.3436 \\
 			(15,20) &       & 0.9160 & 0.2595 & 0.9700 & 0.3070 & 0.9420 & 0.2652 & 0.9780 & 0.3137 \\
 			(20,20) &       & 0.9000 & 0.2463 & 0.9500 & 0.2919 & 0.9060 & 0.2538 & 0.9560 & 0.3008 \\
 			\midrule
 			(10,10) & \multirow{6}[2]{*}{(2,6)} & 0.9480 & 0.2781 & 0.9820 & 0.3314 & 0.9680 & 0.2958 & 0.9780 & 0.3513 \\
 			(10,15) &       & 0.9460 & 0.2402 & 0.9900 & 0.2865 & 0.9700 & 0.2450 & 0.9940 & 0.2923 \\
 			(10,20) &       & 0.9360 & 0.2196 & 0.9760 & 0.2617 & 0.9480 & 0.2179 & 0.9840 & 0.2601 \\
 			(15,15) &       & 0.9260 & 0.2292 & 0.9720 & 0.2735 & 0.9380 & 0.2406 & 0.9760 & 0.2867 \\
 			(15,20) &       & 0.9320 & 0.2084 & 0.9780 & 0.2486 & 0.9500 & 0.2136 & 0.9780 & 0.2549 \\
 			(20,20) &       & 0.9040 & 0.1992 & 0.9540 & 0.2379 & 0.9200 & 0.2077 & 0.9620 & 0.2482 \\
 			\midrule
 			(10,10) & \multirow{6}[2]{*}{(3,4)} & 0.9040 & 0.3687 & 0.9660 & 0.4332 & 0.9100 & 0.3771 & 0.9640 & 0.4427 \\
 			(10,15) &       & 0.9220 & 0.3404 & 0.9720 & 0.4010 & 0.9480 & 0.3479 & 0.9820 & 0.4094 \\
 			(10,20) &       & 0.9280 & 0.3244 & 0.9720 & 0.3826 & 0.9500 & 0.3301 & 0.9860 & 0.3888 \\
 			(15,15) &       & 0.8900 & 0.3136 & 0.9500 & 0.3703 & 0.8960 & 0.3187 & 0.9580 & 0.3761 \\
 			(15,20) &       & 0.8940 & 0.2951 & 0.9640 & 0.3491 & 0.9320 & 0.3000 & 0.9720 & 0.3541 \\
 			(20,20) &       & 0.8760 & 0.2775 & 0.9400 & 0.3283 & 0.8780 & 0.2812 & 0.9500 & 0.3327 \\
 			\midrule
 			(10,10) & \multirow{6}[2]{*}{(3,6)} & 0.9280 & 0.3604 & 0.9680 & 0.4229 & 0.9360 & 0.3763 & 0.9700 & 0.4411 \\
 			(10,15) &       & 0.9400 & 0.3232 & 0.9800 & 0.3796 & 0.9520 & 0.3305 & 0.9860 & 0.3880 \\
 			(10,20) &       & 0.9200 & 0.3018 & 0.9660 & 0.3550 & 0.9440 & 0.3039 & 0.9840 & 0.3573 \\
 			(15,15) &       & 0.8980 & 0.3034 & 0.9580 & 0.3575 & 0.9100 & 0.3137 & 0.9700 & 0.3691 \\
 			(15,20) &       & 0.9120 & 0.2811 & 0.9700 & 0.3316 & 0.9380 & 0.2873 & 0.9780 & 0.3388 \\
 			(20,20) &       & 0.8960 & 0.2668 & 0.9500 & 0.3154 & 0.9020 & 0.2745 & 0.9540 & 0.3245 \\
 			\bottomrule
 		\end{tabular}%
 		\label{unknown.hpd.alpha2.4}%
 	\end{threeparttable}
 \end{table}%

\section*{Case II when $\theta$ is known}
\begin{figure}[htbp]
	\centering
	\includegraphics[width=0.9\linewidth]{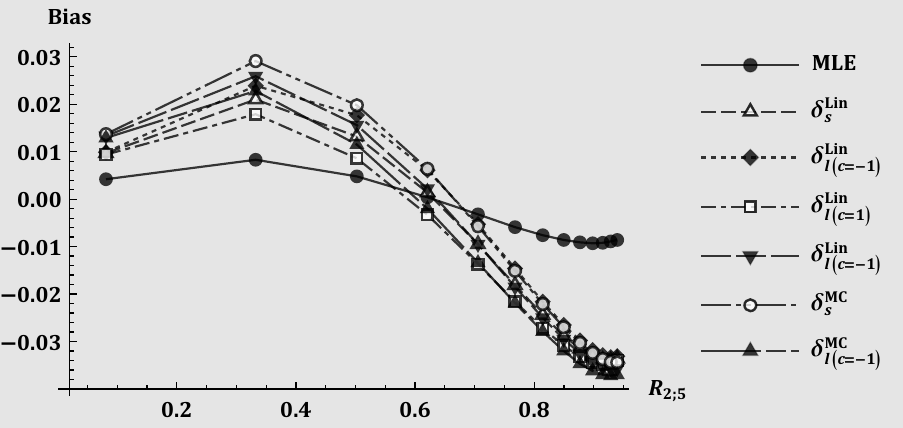}
	\caption{Bias plot of estimators when $\theta$ is known $(\theta=1)$ for $R_{2;5}$}
	\label{fig:bias1knownonlybayes}
\end{figure}

\begin{figure}[htbp]
	\centering
	\includegraphics[width=0.9\linewidth]{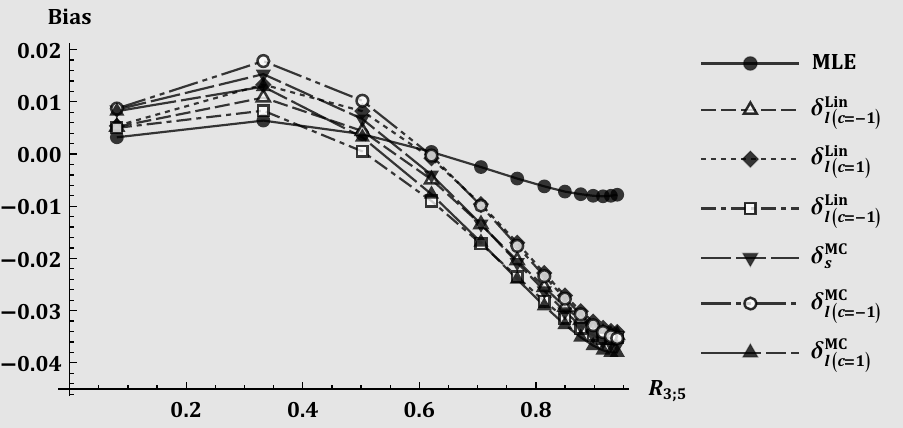}
	\caption{Bias plot of estimators when $\theta$ is known $(\theta=1)$ for $R_{3;5}$}
	\label{fig:bias4knownonlybayes}
\end{figure}

In this case, we again consider distinct components of $(s,k)$ as $\{(1,3),(2,3),(2,5),(3,5), (4,5),$ $(2,6),$ $(3,6),(4,6)\}$  for the numerical aspects. The values of \rsk for given set of $(s,k)$ are  $R_{1;3}=0.90$, $R_{2;3}=0.70$, $R_{2;5}=0.8571$, $R_{3;5}=0.7143$, $R_{4;5}=0.5238$, $R_{2;6}=0.8929$, $R_{3;6}=0.7857$, $R_{4;6}=0.6429$ for $(\alpha_1,\alpha_2,\theta)=(1,2,1).$ All results reported in Table [\ref{known.A.and.mse.alpha1.0.5.prior1}-\ref{known.coverage.proabability}] are for known  $\theta.$ AEs and MSEs of the derived estimators are reported in Table[\ref{known.A.and.mse.alpha1.0.5.prior1}-\ref{known.A.and.mse.alpha1.0.5.prior2.continued}] with different sets of unknown quantities. From these tables we observe that, MSEs of the estimates are decreasing as expected for large number of samples. We also observe that the performance of UMVUE is better with respect to other estimators, in terms of MSEs. Also, We see that Bayes estimators showing lesser MSEs under LINEX loss function. The MCMC Bayes estimator is performing better than Lindley Bayes estimator under SEL function whereas for $c>0$ in LINEX loss function, Lindley Bayes estimator is showing promising results in terms of MSEs. Such type of behavior of Bayes estimator can be seen in Table[\ref{known.A.and.mse.alpha1.0.5.prior1}-\ref{known.A.and.mse.alpha1.0.5.prior2.continued}]. Having said that, UMVUE is still a better estimator with least MSE for MSS reliability and is recommended in practical use. An important observation we make from these tables is that, estimators are showing improvement whenever $m>n$. The Table [\ref{known.coverage.proabability}] reports CPs and ALs of intervals calculated by using the method discussed in Sections [\ref{asymptotic}] and [\ref{bootstrapping}] at $95\%$ level of significance. We see that, the ALs of all kinds of intervals is decreasing as sample sizes increases. Also, for increasing sample sizes the CPs are tending towards the desired level of significance which in this case is $95\%.$ The average Biases of the estimators are exhibited in Figure [\ref{fig:bias1knownonlybayes}-\ref{fig:bias4knownonlybayes}]. The biases are again calculated with respective to the MSS reliability ranging from 0.1 to 0.9 as depicted in the Plots. Figure [\ref{fig:bias1knownonlybayes}] shows the behavior of bias for increasing reliability at  $n=20$ and $m=20$ and $\theta=1.$
Figure [\ref{fig:bias4knownonlybayes}] shows the behavior of bias for increasing reliability at  $n=20$ and $m=25$ and $\theta=1.$
 From Figure [\ref{fig:bias1knownonlybayes}-\ref{fig:bias4knownonlybayes}], it is observed that Bayes estimators are having positive as well as negative average bias for increasing reliability. Also, we see that bias of MLE is remaining close to 0.
\begin{table}[htbp]
	\centering
	\caption{For $(\alpha_1,\alpha_2,\theta)=(0.5,2,1)$ and Prior: $(a_1,a_2)=(2,2)$,$(b_1,b_2)=(3,3).$}
	\begin{tabular}{cccccccccccc}
		\toprule
		\multirow{4}[7]{*}{$(n,m)$} & \multirow{4}[7]{*}{$(s,k)$} & \multicolumn{2}{c}{\multirow{2}[4]{*}{$\hat{R}_{s;k}$}} & \multicolumn{2}{c}{\multirow{2}[4]{*}{$R_{s;k}^{UM}$}} & \multicolumn{6}{c}{ Bayes Lindley Method} \\
		\cmidrule{7-12}          &       & \multicolumn{2}{c}{} & \multicolumn{2}{c}{} & \multicolumn{2}{c}{SEL} & \multicolumn{4}{c}{LINEX} \\
		\cmidrule{3-12}          &       & \multirow{2}[3]{*}{AE} & \multirow{2}[3]{*}{MSE} & \multicolumn{1}{c}{\multirow{2}[3]{*}{AE}} & \multicolumn{1}{c}{\multirow{2}[3]{*}{MSE}} & \multirow{2}[3]{*}{AE} & \multirow{2}[3]{*}{MSE} & \multicolumn{2}{c}{c=-1} & \multicolumn{2}{c}{c=1} \\
		\cmidrule{9-12}          &       &       &       &       &       &       &       & AE    & MSE   & AE    & MSE \\ \midrule
		(10,10) & \multirow{6}[1]{*}{(2,5)} & 0.9461 & 0.0240 & 0.9635 & 0.0016 & 0.9024 & 0.0070 & 0.9020 & 0.0070 & 0.9028 & 0.0070 \\
		(10,15) &       & 0.9472 & 0.0238 & 0.9635 & 0.0012 & 0.9130 & 0.0051 & 0.9134 & 0.0050 & 0.9127 & 0.0053 \\
		(10,20) &       & 0.9475 & 0.0237 & 0.9641 & 0.0010 & 0.9190 & 0.0040 & 0.9196 & 0.0039 & 0.9185 & 0.0041 \\
		(15,15) &       & 0.9506 & 0.0242 & 0.9606 & 0.0013 & 0.9185 & 0.0039 & 0.9186 & 0.0038 & 0.9185 & 0.0039 \\
		(15,20) &       & 0.9510 & 0.0241 & 0.9606 & 0.0012 & 0.9286 & 0.0025 & 0.9288 & 0.0024 & 0.9284 & 0.0025 \\
		(20,20) &       & 0.9529 & 0.0244 & 0.9605 & 0.0009 & 0.9313 & 0.0021 & 0.9314 & 0.0021 & 0.9312 & 0.0022 \\
		\midrule
		(10,10) & \multirow{6}[2]{*}{(3,5)} & 0.8656 & 0.0118 & 0.8836 & 0.0084 & 0.7802 & 0.0177 & 0.7786 & 0.0178 & 0.7817 & 0.0177 \\
		(10,15) &       & 0.8656 & 0.0106 & 0.8843 & 0.0065 & 0.8001 & 0.0128 & 0.8012 & 0.0123 & 0.7992 & 0.0133 \\
		(10,20) &       & 0.8656 & 0.0102 & 0.8848 & 0.0058 & 0.8204 & 0.0084 & 0.8219 & 0.0080 & 0.8191 & 0.0088 \\
		(15,15) &       & 0.8703 & 0.0099 & 0.8822 & 0.0050 & 0.8074 & 0.0105 & 0.8077 & 0.0102 & 0.8073 & 0.0107 \\
		(15,20) &       & 0.8700 & 0.0093 & 0.8822 & 0.0043 & 0.8247 & 0.0076 & 0.8255 & 0.0073 & 0.8240 & 0.0078 \\
		(20,20) &       & 0.8726 & 0.0089 & 0.8816 & 0.0037 & 0.8287 & 0.0065 & 0.8292 & 0.0063 & 0.8283 & 0.0066 \\
		\midrule
		(10,10) & \multirow{6}[2]{*}{(4,5)} & 0.7155 & 0.0206 & 0.7294 & 0.0138 & 0.5950 & 0.0244 & 0.5910 & 0.0249 & 0.5984 & 0.0241 \\
		(10,15) &       & 0.7130 & 0.0189 & 0.7296 & 0.0114 & 0.6294 & 0.0164 & 0.6310 & 0.0156 & 0.6282 & 0.0171 \\
		(10,20) &       & 0.7122 & 0.0183 & 0.7304 & 0.0101 & 0.6473 & 0.0124 & 0.6502 & 0.0116 & 0.6447 & 0.0131 \\
		(15,15) &       & 0.7175 & 0.0161 & 0.7246 & 0.0101 & 0.6357 & 0.0144 & 0.6360 & 0.0139 & 0.6356 & 0.0147 \\
		(15,20) &       & 0.7159 & 0.0153 & 0.7242 & 0.0089 & 0.6567 & 0.0108 & 0.6584 & 0.0104 & 0.6553 & 0.0113 \\
		(20,20) &       & 0.7182 & 0.0137 & 0.7236 & 0.0074 & 0.6578 & 0.0092 & 0.6590 & 0.0089 & 0.6569 & 0.0095 \\
		\midrule
		(10,10) & \multirow{6}[2]{*}{(2,6)} & 0.9634 & 0.0284 & 0.9747 & 0.0037 & 0.9300 & 0.0050 & 0.9298 & 0.0049 & 0.9301 & 0.0050 \\
		(10,15) &       & 0.9645 & 0.0284 & 0.9736 & 0.0028 & 0.9371 & 0.0036 & 0.9373 & 0.0036 & 0.9368 & 0.0037 \\
		(10,20) &       & 0.9649 & 0.0285 & 0.9718 & 0.0039 & 0.9413 & 0.0031 & 0.9417 & 0.0029 & 0.9409 & 0.0032 \\
		(15,15) &       & 0.9675 & 0.0290 & 0.9768 & 0.0006 & 0.9437 & 0.0027 & 0.9438 & 0.0026 & 0.9437 & 0.0027 \\
		(15,20) &       & 0.9680 & 0.0290 & 0.9770 & 0.0005 & 0.9487 & 0.0019 & 0.9489 & 0.0019 & 0.9486 & 0.0020 \\
		(20,20) &       & 0.9696 & 0.0294 & 0.9762 & 0.0005 & 0.9518 & 0.0016 & 0.9519 & 0.0016 & 0.9517 & 0.0017 \\
		\midrule
		(10,10) & \multirow{6}[2]{*}{(3,6)} & 0.9117 & 0.0176 & 0.9371 & 0.0030 & 0.8381 & 0.0157 & 0.8376 & 0.0155 & 0.8385 & 0.0159 \\
		(10,15) &       & 0.9125 & 0.0169 & 0.9389 & 0.0022 & 0.8631 & 0.0093 & 0.8638 & 0.0089 & 0.8626 & 0.0096 \\
		(10,20) &       & 0.9127 & 0.0167 & 0.9399 & 0.0019 & 0.8724 & 0.0072 & 0.8735 & 0.0069 & 0.8714 & 0.0076 \\
		(15,15) &       & 0.9171 & 0.0165 & 0.9290 & 0.0035 & 0.8719 & 0.0070 & 0.8720 & 0.0068 & 0.8718 & 0.0071 \\
		(15,20) &       & 0.9195 & 0.0166 & 0.9291 & 0.0029 & 0.8790 & 0.0056 & 0.8796 & 0.0055 & 0.8785 & 0.0058 \\
		(20,20) &       & 0.8195 & 0.0109 & 0.9289 & 0.0022 & 0.8873 & 0.0041 & 0.8875 & 0.0040 & 0.8870 & 0.0042 \\
		\midrule
		(10,10) & \multirow{6}[2]{*}{(4,6)} & 0.8187 & 0.0093 & 0.8383 & 0.0120 & 0.7094 & 0.0253 & 0.7072 & 0.0254 & 0.7115 & 0.0252 \\
		(10,15) &       & 0.8184 & 0.0087 & 0.8356 & 0.0120 & 0.7456 & 0.0157 & 0.7469 & 0.0150 & 0.7446 & 0.0163 \\
		(10,20) &       & 0.8237 & 0.0077 & 0.8347 & 0.0123 & 0.7619 & 0.0116 & 0.7641 & 0.0110 & 0.7599 & 0.0123 \\
		(15,15) &       & 0.8230 & 0.0068 & 0.8360 & 0.0073 & 0.7496 & 0.0135 & 0.7499 & 0.0132 & 0.7495 & 0.0138 \\
		(15,20) &       & 0.8256 & 0.0060 & 0.8361 & 0.0062 & 0.7688 & 0.0095 & 0.7701 & 0.0092 & 0.7678 & 0.0099 \\
		(20,20) &       & 0.8256 & 0.0687 & 0.8342 & 0.0055 & 0.7729 & 0.0086 & 0.7737 & 0.0084 & 0.7723 & 0.0089 \\
		\bottomrule
	\end{tabular}%
	
	\label{known.A.and.mse.alpha1.0.5.prior1}%
\end{table}%

\begin{table}[htbp]
	\centering
	\caption{\textit{continued}}
	\begin{tabular}{ccccccccccc}
		\toprule
		\multicolumn{2}{c}{Bayes Lindley} &       & \multicolumn{8}{c}{Bayes MCMC } \\
		\midrule
		\multicolumn{2}{c}{LINEX} &       & \multicolumn{2}{c}{SEL} & \multicolumn{6}{c}{LINEX} \\
		\midrule
		\multicolumn{2}{c}{c=1.5} &       & \multirow{2}[4]{*}{AE} & \multirow{2}[4]{*}{MSE} & \multicolumn{2}{c}{c=-1} & \multicolumn{2}{c}{c=1} & \multicolumn{2}{c}{c=1.5} \\
		\cmidrule{1-2}\cmidrule{6-11}    AE    & MSE   &       &       &       & AE    & MSE   & AE    & MSE   & AE    & MSE \\
		\midrule
		0.9030 & 0.0070 &       & 0.8974 & 0.0068 & 0.8998 & 0.0063 & 0.8949 & 0.0073 & 0.8936 & 0.0075 \\
		0.9126 & 0.0053 &       & 0.9126 & 0.0045 & 0.9143 & 0.0043 & 0.9109 & 0.0048 & 0.9100 & 0.0050 \\
		0.9182 & 0.0042 &       & 0.9208 & 0.0035 & 0.9222 & 0.0033 & 0.9195 & 0.0037 & 0.9188 & 0.0038 \\
		0.9185 & 0.0039 &       & 0.9159 & 0.0039 & 0.9173 & 0.0037 & 0.9145 & 0.0041 & 0.9138 & 0.0042 \\
		0.9283 & 0.0025 &       & 0.9238 & 0.0029 & 0.9249 & 0.0028 & 0.9228 & 0.0030 & 0.9222 & 0.0031 \\
		0.9312 & 0.0022 &       & 0.9261 & 0.0026 & 0.9270 & 0.0025 & 0.9252 & 0.0027 & 0.9247 & 0.0027 \\
		\cmidrule{1-2}\cmidrule{4-11}    0.7824 & 0.0176 &       & 0.7899 & 0.0144 & 0.7945 & 0.0133 & 0.7851 & 0.0156 & 0.7826 & 0.0162 \\
		0.7988 & 0.0135 &       & 0.8114 & 0.0102 & 0.8149 & 0.0096 & 0.8077 & 0.0110 & 0.8058 & 0.0114 \\
		0.8185 & 0.0090 &       & 0.8235 & 0.0082 & 0.8265 & 0.0077 & 0.8204 & 0.0088 & 0.8188 & 0.0091 \\
		0.8072 & 0.0107 &       & 0.8145 & 0.0091 & 0.8175 & 0.0086 & 0.8114 & 0.0097 & 0.8099 & 0.0100 \\
		0.8237 & 0.0079 &       & 0.8260 & 0.0071 & 0.8285 & 0.0067 & 0.8235 & 0.0075 & 0.8222 & 0.0077 \\
		0.8281 & 0.0067 &       & 0.8285 & 0.0064 & 0.8306 & 0.0061 & 0.8262 & 0.0067 & 0.8251 & 0.0069 \\
		\cmidrule{1-2}\cmidrule{4-11}    0.6000 & 0.0240 &       & 0.6237 & 0.0182 & 0.6298 & 0.0169 & 0.6176 & 0.0195 & 0.6144 & 0.0202 \\
		0.6277 & 0.0174 &       & 0.6477 & 0.0136 & 0.6528 & 0.0127 & 0.6426 & 0.0145 & 0.6400 & 0.0150 \\
		0.6436 & 0.0135 &       & 0.6617 & 0.0113 & 0.6662 & 0.0107 & 0.6572 & 0.0120 & 0.6549 & 0.0124 \\
		0.6356 & 0.0149 &       & 0.6490 & 0.0123 & 0.6533 & 0.0116 & 0.6447 & 0.0130 & 0.6425 & 0.0134 \\
		0.6547 & 0.0115 &       & 0.6621 & 0.0100 & 0.6658 & 0.0095 & 0.6583 & 0.0105 & 0.6564 & 0.0108 \\
		0.6565 & 0.0097 &       & 0.6636 & 0.0091 & 0.6669 & 0.0087 & 0.6603 & 0.0096 & 0.6586 & 0.0098 \\
		\cmidrule{1-2}\cmidrule{4-11}    0.9302 & 0.0050 &       & 0.9228 & 0.0049 & 0.9246 & 0.0046 & 0.9209 & 0.0052 & 0.9198 & 0.0054 \\
		0.9367 & 0.0038 &       & 0.9359 & 0.0032 & 0.9372 & 0.0030 & 0.9347 & 0.0034 & 0.9340 & 0.0035 \\
		0.9407 & 0.0032 &       & 0.9430 & 0.0024 & 0.9439 & 0.0023 & 0.9420 & 0.0026 & 0.9415 & 0.0026 \\
		0.9436 & 0.0027 &       & 0.9392 & 0.0027 & 0.9401 & 0.0026 & 0.9382 & 0.0029 & 0.9376 & 0.0029 \\
		0.9485 & 0.0020 &       & 0.9459 & 0.0020 & 0.9466 & 0.0019 & 0.9451 & 0.0021 & 0.9447 & 0.0021 \\
		0.9517 & 0.0017 &       & 0.9480 & 0.0017 & 0.9486 & 0.0017 & 0.9473 & 0.0018 & 0.9470 & 0.0018 \\
		\cmidrule{1-2}\cmidrule{4-11}    0.8388 & 0.0160 &       & 0.8466 & 0.0115 & 0.8503 & 0.0106 & 0.8427 & 0.0124 & 0.8406 & 0.0129 \\
		0.8623 & 0.0097 &       & 0.8659 & 0.0079 & 0.8687 & 0.0074 & 0.8630 & 0.0085 & 0.8615 & 0.0088 \\
		0.8709 & 0.0077 &       & 0.8766 & 0.0063 & 0.8789 & 0.0059 & 0.8743 & 0.0067 & 0.8730 & 0.0069 \\
		0.8530 & 0.0110 &       & 0.8498 & 0.0101 & 0.8531 & 0.0094 & 0.8464 & 0.0108 & 0.8447 & 0.0112 \\
		0.8718 & 0.0071 &       & 0.8695 & 0.0070 & 0.8718 & 0.0066 & 0.8671 & 0.0074 & 0.8659 & 0.0076 \\
		0.8783 & 0.0059 &       & 0.8797 & 0.0053 & 0.8816 & 0.0050 & 0.8778 & 0.0056 & 0.8768 & 0.0058 \\
		0.8869 & 0.0042 &       & 0.8823 & 0.0048 & 0.8839 & 0.0045 & 0.8806 & 0.0050 & 0.8798 & 0.0051 \\
		\cmidrule{1-2}\cmidrule{4-11}    0.7124 & 0.0252 &       & 0.7332 & 0.0177 & 0.7389 & 0.0164 & 0.7274 & 0.0192 & 0.7244 & 0.0200 \\
		0.7442 & 0.0166 &       & 0.7569 & 0.0129 & 0.7614 & 0.0120 & 0.7522 & 0.0138 & 0.7499 & 0.0143 \\
		0.7590 & 0.0126 &       & 0.7704 & 0.0105 & 0.7744 & 0.0098 & 0.7664 & 0.0112 & 0.7644 & 0.0116 \\
		0.7494 & 0.0139 &       & 0.7596 & 0.0116 & 0.7634 & 0.0108 & 0.7556 & 0.0123 & 0.7536 & 0.0127 \\
		0.7673 & 0.0100 &       & 0.7724 & 0.0092 & 0.7756 & 0.0086 & 0.7690 & 0.0097 & 0.7674 & 0.0100 \\
		0.7720 & 0.0090 &       & 0.7746 & 0.0083 & 0.7775 & 0.0079 & 0.7717 & 0.0088 & 0.7702 & 0.0090 \\
		\bottomrule    \end{tabular}%
	\label{known.A.and.mse.alpha1.0.5.prior1.conti}%
\end{table}%
\begin{table}[htbp]
	\centering
	\caption{For $(\alpha_1,\alpha_2,\theta)=(0.5,2,1)$ and Prior: $(a_1,a_2)=(2,2)$,$(b_1,b_2)=(4,4).$}
	\begin{tabular}{cccccccccc}
		\toprule
		\multirow{4}[8]{*}{$(n,m)$} & \multirow{4}[8]{*}{$(s,k)$} & \multicolumn{8}{c}{ Bayes Lindley Method} \\
		\cmidrule{3-10}          &       & \multicolumn{2}{c}{SEL} & \multicolumn{6}{c}{LINEX} \\
		\cmidrule{3-10}          &       & \multirow{2}[4]{*}{AE} & \multirow{2}[4]{*}{MSE} & \multicolumn{2}{c}{c=-1} & \multicolumn{2}{c}{c=1} & \multicolumn{2}{c}{c=1.5} \\
		\cmidrule{5-10}          &       &       &       & AE    & MSE   & AE    & MSE   & AE    & MSE \\
		\midrule
		(10,10) & \multirow{6}[2]{*}{(2,5)} & 0.8813 & 0.0113 & 0.8802 & 0.0114 & 0.8823 & 0.0112 & 0.8828 & 0.0111 \\
		(10,15) &       & 0.9085 & 0.0053 & 0.9085 & 0.0052 & 0.9085 & 0.0054 & 0.9085 & 0.0054 \\
		(10,20) &       & 0.9151 & 0.0042 & 0.9155 & 0.0041 & 0.9147 & 0.0043 & 0.9145 & 0.0044 \\
		(15,15) &       & 0.9097 & 0.0048 & 0.9094 & 0.0048 & 0.9100 & 0.0048 & 0.9102 & 0.0048 \\
		(15,20) &       & 0.9191 & 0.0035 & 0.9191 & 0.0034 & 0.9190 & 0.0035 & 0.9190 & 0.0036 \\
		(20,20) &       & 0.9267 & 0.0028 & 0.9266 & 0.0027 & 0.9268 & 0.0028 & 0.9268 & 0.0028 \\
		\midrule
		(10,10) & \multirow{6}[2]{*}{(3,5)} & 0.7556 & 0.0240 & 0.7504 & 0.0254 & 0.7602 & 0.0230 & 0.7623 & 0.0225 \\
		(10,15) &       & 0.7919 & 0.0139 & 0.7917 & 0.0136 & 0.7922 & 0.0141 & 0.7924 & 0.0142 \\
		(10,20) &       & 0.8080 & 0.0104 & 0.8092 & 0.0100 & 0.8069 & 0.0109 & 0.8064 & 0.0111 \\
		(15,15) &       & 0.7975 & 0.0122 & 0.7962 & 0.0123 & 0.7988 & 0.0121 & 0.7994 & 0.0120 \\
		(15,20) &       & 0.8109 & 0.0092 & 0.8111 & 0.0091 & 0.8108 & 0.0094 & 0.8108 & 0.0095 \\
		(20,20) &       & 0.8164 & 0.0078 & 0.8162 & 0.0078 & 0.8167 & 0.0079 & 0.8169 & 0.0079 \\
		\midrule
		(10,10) & \multirow{6}[2]{*}{(4,5)} & 0.5521 & 0.0374 & 0.5391 & 0.0418 & 0.5621 & 0.0346 & 0.5663 & 0.0335 \\
		(10,15) &       & 0.6079 & 0.0205 & 0.6070 & 0.0202 & 0.6089 & 0.0209 & 0.6093 & 0.0210 \\
		(10,20) &       & 0.6352 & 0.0139 & 0.6372 & 0.0132 & 0.6336 & 0.0146 & 0.6329 & 0.0149 \\
		(15,15) &       & 0.6083 & 0.0193 & 0.6052 & 0.0196 & 0.6110 & 0.0191 & 0.6123 & 0.0190 \\
		(15,20) &       & 0.6368 & 0.0132 & 0.6370 & 0.0129 & 0.6367 & 0.0135 & 0.6367 & 0.0136 \\
		(20,20) &       & 0.6367 & 0.0122 & 0.6361 & 0.0121 & 0.6373 & 0.0123 & 0.6377 & 0.0124 \\
		\midrule
		(10,10) & \multirow{6}[2]{*}{(2,6)} & 0.9171 & 0.0071 & 0.9164 & 0.0072 & 0.9177 & 0.0070 & 0.9180 & 0.0070 \\
		(10,15) &       & 0.9365 & 0.0035 & 0.9365 & 0.0034 & 0.9364 & 0.0035 & 0.9364 & 0.0036 \\
		(10,20) &       & 0.9420 & 0.0028 & 0.9423 & 0.0027 & 0.9417 & 0.0028 & 0.9416 & 0.0029 \\
		(15,15) &       & 0.9364 & 0.0032 & 0.9363 & 0.0032 & 0.9366 & 0.0032 & 0.9367 & 0.0032 \\
		(15,20) &       & 0.9448 & 0.0022 & 0.9448 & 0.0021 & 0.9448 & 0.0022 & 0.9448 & 0.0022 \\
		(20,20) &       & 0.9468 & 0.0019 & 0.9468 & 0.0019 & 0.9468 & 0.0019 & 0.9469 & 0.0019 \\
		\midrule
		(10,10) & \multirow{6}[2]{*}{(3,6)} & 0.8210 & 0.0188 & 0.8183 & 0.0193 & 0.8234 & 0.0183 & 0.8245 & 0.0181 \\
		(10,15) &       & 0.8549 & 0.0101 & 0.8548 & 0.0099 & 0.8550 & 0.0103 & 0.8551 & 0.0103 \\
		(10,20) &       & 0.8683 & 0.0075 & 0.8692 & 0.0072 & 0.8676 & 0.0078 & 0.8672 & 0.0079 \\
		(15,15) &       & 0.8567 & 0.0093 & 0.8559 & 0.0093 & 0.8575 & 0.0092 & 0.8578 & 0.0092 \\
		(15,20) &       & 0.8720 & 0.0066 & 0.8722 & 0.0065 & 0.8719 & 0.0067 & 0.8719 & 0.0068 \\
		(20,20) &       & 0.8741 & 0.0059 & 0.8739 & 0.0058 & 0.8743 & 0.0059 & 0.8744 & 0.0059 \\
		\midrule
		(10,10) & \multirow{6}[2]{*}{(4,6)} & 0.6889 & 0.0313 & 0.6816 & 0.0333 & 0.6951 & 0.0298 & 0.6979 & 0.0292 \\
		(10,15) &       & 0.7271 & 0.0186 & 0.7266 & 0.0182 & 0.7276 & 0.0189 & 0.7279 & 0.0190 \\
		(10,20) &       & 0.7498 & 0.0133 & 0.7515 & 0.0126 & 0.7484 & 0.0139 & 0.7478 & 0.0142 \\
		(15,15) &       & 0.7313 & 0.0164 & 0.7292 & 0.0166 & 0.7334 & 0.0162 & 0.7344 & 0.0160 \\
		(15,20) &       & 0.7585 & 0.0110 & 0.7586 & 0.0108 & 0.7585 & 0.0112 & 0.7586 & 0.0113 \\
		(20,20) &       & 0.7539 & 0.0116 & 0.7536 & 0.0115 & 0.7543 & 0.0117 & 0.7545 & 0.0117 \\
		\bottomrule
	\end{tabular}%
	\label{known.A.and.mse.alpha1.0.5.prior2}%
\end{table}%
\begin{table}[htbp]
	\centering
	\caption{\textit{continued}}
	\begin{tabular}{cccccccccc}
		\toprule
		\multirow{4}[7]{*}{$(n,m)$} & \multirow{4}[7]{*}{$(s,k)$} & \multicolumn{8}{c}{ Bayes MCMC Method} \\
		\cmidrule{3-10}          &       & \multicolumn{2}{c}{SEL} & \multicolumn{6}{c}{LINEX} \\
		\cmidrule{3-10}          &       & \multirow{2}[3]{*}{AE} & \multirow{2}[3]{*}{MSE} & \multicolumn{2}{c}{c=-1} & \multicolumn{2}{c}{c=1} & \multicolumn{2}{c}{c=1.5} \\
		\cmidrule{5-10}          &       &       &       & AE    & MSE   & AE    & MSE   & AE    & MSE \\ \midrule
		(10,10) & \multirow{6}[2]{*}{(2,5)} & 0.8863 & 0.0083 & 0.8890 & 0.0078 & 0.8835 & 0.0089 & 0.8820 & 0.0092 \\
		(10,15) &       & 0.9066 & 0.0051 & 0.9084 & 0.0048 & 0.9047 & 0.0054 & 0.9037 & 0.0056 \\
		(10,20) &       & 0.9175 & 0.0037 & 0.9189 & 0.0035 & 0.9161 & 0.0039 & 0.9154 & 0.0040 \\
		(15,15) &       & 0.9078 & 0.0048 & 0.9093 & 0.0045 & 0.9062 & 0.0050 & 0.9054 & 0.0051 \\
		(15,20) &       & 0.9185 & 0.0033 & 0.9196 & 0.0032 & 0.9173 & 0.0035 & 0.9167 & 0.0036 \\
		(20,20) &       & 0.9197 & 0.0031 & 0.9207 & 0.0030 & 0.9187 & 0.0032 & 0.9182 & 0.0033 \\
		\midrule
		(10,10) & \multirow{6}[2]{*}{(3,5)} & 0.7730 & 0.0175 & 0.7780 & 0.0162 & 0.7678 & 0.0189 & 0.7651 & 0.0196 \\
		(10,15) &       & 0.8014 & 0.0114 & 0.8052 & 0.0107 & 0.7975 & 0.0123 & 0.7955 & 0.0127 \\
		(10,20) &       & 0.8176 & 0.0086 & 0.8208 & 0.0081 & 0.8144 & 0.0092 & 0.8128 & 0.0095 \\
		(15,15) &       & 0.8012 & 0.0109 & 0.8044 & 0.0103 & 0.7979 & 0.0116 & 0.7962 & 0.0120 \\
		(15,20) &       & 0.8170 & 0.0081 & 0.8196 & 0.0076 & 0.8143 & 0.0085 & 0.8129 & 0.0088 \\
		(20,20) &       & 0.8176 & 0.0076 & 0.8200 & 0.0072 & 0.8153 & 0.0080 & 0.8141 & 0.0082 \\
		\midrule
		(10,10) & \multirow{6}[2]{*}{(4,5)} & 0.6033 & 0.0218 & 0.6096 & 0.0202 & 0.5970 & 0.0234 & 0.5938 & 0.0242 \\
		(10,15) &       & 0.6349 & 0.0149 & 0.6401 & 0.0139 & 0.6297 & 0.0160 & 0.6270 & 0.0165 \\
		(10,20) &       & 0.6537 & 0.0117 & 0.6583 & 0.0110 & 0.6491 & 0.0125 & 0.6467 & 0.0129 \\
		(15,15) &       & 0.6323 & 0.0145 & 0.6367 & 0.0136 & 0.6278 & 0.0154 & 0.6256 & 0.0158 \\
		(15,20) &       & 0.6503 & 0.0111 & 0.6541 & 0.0105 & 0.6464 & 0.0118 & 0.6445 & 0.0121 \\
		(20,20) &       & 0.6497 & 0.0106 & 0.6531 & 0.0101 & 0.6463 & 0.0112 & 0.6446 & 0.0115 \\
		\midrule
		(10,10) & \multirow{6}[2]{*}{(2,6)} & 0.9135 & 0.0060 & 0.9156 & 0.0057 & 0.9113 & 0.0065 & 0.9101 & 0.0067 \\
		(10,15) &       & 0.9311 & 0.0036 & 0.9324 & 0.0034 & 0.9297 & 0.0038 & 0.9289 & 0.0039 \\
		(10,20) &       & 0.9404 & 0.0025 & 0.9414 & 0.0024 & 0.9393 & 0.0027 & 0.9388 & 0.0028 \\
		(15,15) &       & 0.9325 & 0.0033 & 0.9336 & 0.0032 & 0.9314 & 0.0035 & 0.9308 & 0.0036 \\
		(15,20) &       & 0.9416 & 0.0023 & 0.9425 & 0.0022 & 0.9408 & 0.0024 & 0.9404 & 0.0024 \\
		(20,20) &       & 0.9429 & 0.0021 & 0.9436 & 0.0020 & 0.9422 & 0.0022 & 0.9418 & 0.0022 \\
		\midrule
		(10,10) & \multirow{6}[2]{*}{(3,6)} & 0.8319 & 0.0140 & 0.8360 & 0.0130 & 0.8276 & 0.0151 & 0.8253 & 0.0157 \\
		(10,15) &       & 0.8576 & 0.0089 & 0.8605 & 0.0083 & 0.8545 & 0.0095 & 0.8529 & 0.0099 \\
		(10,20) &       & 0.8718 & 0.0066 & 0.8742 & 0.0062 & 0.8694 & 0.0070 & 0.8681 & 0.0072 \\
		(15,10) &       & 0.8319 & 0.0133 & 0.8355 & 0.0124 & 0.8281 & 0.0142 & 0.8261 & 0.0147 \\
		(15,15) &       & 0.8582 & 0.0084 & 0.8607 & 0.0079 & 0.8556 & 0.0089 & 0.8543 & 0.0092 \\
		(15,20) &       & 0.8722 & 0.0061 & 0.8742 & 0.0058 & 0.8702 & 0.0064 & 0.8691 & 0.0066 \\
		(20,20) &       & 0.8733 & 0.0057 & 0.8751 & 0.0054 & 0.8715 & 0.0060 & 0.8706 & 0.0061 \\
		\midrule
		(10,10) & \multirow{6}[2]{*}{(4,6)} & 0.7140 & 0.0214 & 0.7200 & 0.0198 & 0.7078 & 0.0232 & 0.7046 & 0.0241 \\
		(10,15) &       & 0.7453 & 0.0143 & 0.7500 & 0.0133 & 0.7404 & 0.0154 & 0.7379 & 0.0160 \\
		(10,20) &       & 0.7634 & 0.0110 & 0.7674 & 0.0102 & 0.7593 & 0.0118 & 0.7572 & 0.0122 \\
		(15,15) &       & 0.7441 & 0.0138 & 0.7481 & 0.0129 & 0.7400 & 0.0147 & 0.7379 & 0.0152 \\
		(15,20) &       & 0.7617 & 0.0103 & 0.7651 & 0.0097 & 0.7582 & 0.0110 & 0.7565 & 0.0113 \\
		(20,20) &       & 0.7620 & 0.0098 & 0.7650 & 0.0093 & 0.7589 & 0.0104 & 0.7573 & 0.0107 \\
		\bottomrule
	\end{tabular}%
	\label{known.A.and.mse.alpha1.0.5.prior2.continued}%
\end{table}%

\begin{landscape}
\begin{table}[htbp]
	\centering
	\caption{For $(\alpha_1,\alpha_2,\theta)=(0.5,2,1)$, $Prior_1$: $(a_1,a_2)=(1,1)$,$(b_1,b_2)=(1.5,1.5)$ and $Prior_2$: $(a_1,a_2)=(1.5,1.5)$,$(b_1,b_2)=(1.5,1.5).$ }
	\begin{tabular}{cccccccccccccc}
		\toprule
		\multirow{3}[6]{*}{$(n,m)$} & \multirow{3}[6]{*}{$(s,k)$} &       &       & \multicolumn{6}{c}{Boot Intervals}            & \multicolumn{4}{c}{HPD Interval} \\
		\cmidrule{3-14}          &       & \multicolumn{2}{c}{Asymptotic} & \multicolumn{2}{c}{Boot normal} & \multicolumn{2}{c}{Boot percentile} & \multicolumn{2}{c}{Boot T} & \multicolumn{2}{c}{$Prior_1$} & \multicolumn{2}{c}{$Prior_2$} \\
		\cmidrule{3-14}          &       & CP    & AL    & CP    & AL    & CP    & AL    & CP    & AL    & CP    & AL    & CP    & AL \\
		\midrule
		(10,10) & \multirow{6}[1]{*}{(1,3)} & 0.8700 & 0.1229 & 0.9340 & 0.1524 & 0.9410 & 0.1436 & 0.9670 & 0.2228 & 0.9720 & 0.1456 & 0.9740 & 0.1423 \\
		(10,15) &       & 0.8980 & 0.1114 & 0.9510 & 0.1411 & 0.9530 & 0.1340 & 0.9870 & 0.2081 & 0.9800 & 0.1210 & 0.9700 & 0.1210 \\
		(10,20) &       & 0.8890 & 0.1036 & 0.9520 & 0.1345 & 0.9460 & 0.1280 & 0.9830 & 0.1983 & 0.9580 & 0.1084 & 0.9620 & 0.1098 \\
		(15,15) &       & 0.8800 & 0.0945 & 0.9260 & 0.1127 & 0.9370 & 0.1077 & 0.9680 & 0.1607 & 0.9660 & 0.1084 & 0.9640 & 0.1067 \\
		(15,20) &       & 0.8770 & 0.0843 & 0.9320 & 0.1059 & 0.9470 & 0.1018 & 0.9770 & 0.1514 & 0.9680 & 0.0969 & 0.9640 & 0.0965 \\
		(20,20) &       & 0.8800 & 0.0800 & 0.9280 & 0.0917 & 0.9490 & 0.0885 & 0.9720 & 0.1279 & 0.9560 & 0.0894 & 0.9560 & 0.0885 \\\midrule
		(10,10) & \multirow{6}[1]{*}{(2,3)} & 0.9140 & 0.3073 & 0.9170 & 0.3196 & 0.9410 & 0.3099 & 0.9540 & 0.4058 & 0.9520 & 0.3171 & 0.9560 & 0.3113 \\
		(10,15) &       & 0.9060 & 0.2888 & 0.9340 & 0.3004 & 0.9530 & 0.2930 & 0.9770 & 0.3878 & 0.9560 & 0.2817 & 0.9520 & 0.2792 \\
		(10,20) &       & 0.9290 & 0.2711 & 0.9360 & 0.2884 & 0.9460 & 0.2816 & 0.9720 & 0.3742 & 0.9460 & 0.2619 & 0.9460 & 0.2614 \\
		(15,15) &       & 0.9140 & 0.2526 & 0.9190 & 0.2593 & 0.9370 & 0.2536 & 0.9580 & 0.3217 & 0.9520 & 0.2573 & 0.9500 & 0.2536 \\
		(15,20) &       & 0.9380 & 0.2399 & 0.9270 & 0.2461 & 0.9470 & 0.2416 & 0.9700 & 0.3076 & 0.9500 & 0.2380 & 0.9500 & 0.2363 \\
		(20,20) &       & 0.9280 & 0.2189 & 0.9240 & 0.2239 & 0.9490 & 0.2198 & 0.9650 & 0.2718 & 0.9400 & 0.2225 & 0.9340 & 0.2202 \\
		\midrule
		(10,10) & \multirow{6}[1]{*}{(2,5)} & 0.8670 & 0.1765 & 0.9280 & 0.2081 & 0.9410 & 0.1963 & 0.9640 & 0.3005 & 0.9680 & 0.2004 & 0.9720 & 0.1961 \\
		(10,15) &       & 0.8700 & 0.1578 & 0.9480 & 0.1943 & 0.9530 & 0.1848 & 0.9840 & 0.2833 & 0.9740 & 0.1688 & 0.9680 & 0.1687 \\
		(10,20) &       & 0.8770 & 0.1495 & 0.9460 & 0.1858 & 0.9460 & 0.1771 & 0.9800 & 0.2713 & 0.9580 & 0.1523 & 0.9600 & 0.1541 \\
		(15,15) &       & 0.8760 & 0.1367 & 0.9240 & 0.1582 & 0.9370 & 0.1511 & 0.9670 & 0.2233 & 0.9640 & 0.1525 & 0.9600 & 0.1502 \\
		(15,20) &       & 0.8940 & 0.1308 & 0.9310 & 0.1493 & 0.9470 & 0.1435 & 0.9760 & 0.2116 & 0.9660 & 0.1374 & 0.9600 & 0.1368 \\
		(20,20) &       & 0.8910 & 0.1180 & 0.9250 & 0.1306 & 0.9490 & 0.1260 & 0.9710 & 0.1808 & 0.9540 & 0.1275 & 0.9540 & 0.1262 \\ \midrule
		(10,10) & \multirow{6}[1]{*}{(3,5)} & 0.8740 & 0.3094 & 0.9130 & 0.3286 & 0.9410 & 0.3160 & 0.9530 & 0.4275 & 0.9560 & 0.3246 & 0.9560 & 0.3187 \\
		(10,15) &       & 0.9030 & 0.2861 & 0.9320 & 0.3098 & 0.9530 & 0.3000 & 0.9770 & 0.4090 & 0.9540 & 0.2871 & 0.9480 & 0.2851 \\
		(10,20) &       & 0.9080 & 0.2816 & 0.9310 & 0.2978 & 0.9460 & 0.2890 & 0.9720 & 0.3949 & 0.9460 & 0.2662 & 0.9460 & 0.2665 \\
		(15,15) &       & 0.9090 & 0.2617 & 0.9140 & 0.2663 & 0.9370 & 0.2587 & 0.9580 & 0.3388 & 0.9500 & 0.2633 & 0.9500 & 0.2596 \\
		(15,20) &       & 0.9010 & 0.2380 & 0.9220 & 0.2532 & 0.9470 & 0.2472 & 0.9700 & 0.3241 & 0.9500 & 0.2431 & 0.9500 & 0.2416 \\
		(20,20) &       & 0.9220 & 0.2227 & 0.9200 & 0.2296 & 0.9490 & 0.2243 & 0.9660 & 0.2860 & 0.9420 & 0.2277 & 0.9380 & 0.2254 \\
		\bottomrule    \end{tabular}%
	\label{known.coverage.proabability}%
\end{table}%
\end{landscape}


\section{An Illustrative Example}
In this section, two real datasets of breakdown times of several specimens of electrical insulating fluid are considered.   \cite{nelson1972graphical} obtained these datasets by conducting an experiment to investigate the distribution of breakdown times of several specimens of a particular type of electrical insulating fluid at constant voltage ranging from 26 kilovolts (kV) to 38 kV (see also \cite{azhad2019common}). The datasets are
\begin{center}
	\begin{tabular}{lcl}
		\toprule
		\multirow{2}{*}{Data I(strength)}  &\multirow{2}{*}{:} & 0.40, 82.85, 9.88, 89.29, 215.10, 2.75, 0.79, 15.93, 3.91,	0.27, 0.69, 100.58 \\
		&   & 27.80, 13.95, 53,24                                                           \\ \midrule
		\multirow{1}{*}{Data II(stress)} & \multirow{1}{*}{:} & 0.47, 0.73, 1.40, 0.74, 0.39, 1.13, 0.09, 2.38\\ \bottomrule
	\end{tabular}
\end{center} 
\noindent
With the help of Kolmogorov-Smirnov test we find that Data I supports $P(0.3,0.8)$ distribution with  KS distance $0.19$ and  $p$-value $0.5711$ and Data II supports $P(1.4,0.8)$ distribution with KS distance $0.3$ and  $p$-value $0.4.$ In this problem, we assume that electrical fluid of specimen considered to be good if 2 out of 4 specimens are functioning properly at constant voltage. Form Data I and Data II, two sets of upper record values $\boldsymbol{R}= (0.40,  82.85,  89.29, 215.10)$ and $\boldsymbol{S}= (0.47, 0.73, 1.40, 2.38)$ are obtained, respectively. From \BR and \BS, we find that $n=4$, $m=4$, $\hat{\theta}=0.4$, $\hat{\alpha_1}=0.64$, $\hat{\alpha2}=2.24$.  The ML estimate of \rsk for the given datasets is $\hat{R}_{2;4}=0.9105$. For prior distribution, we consider $(a_1,a_2,a_3)=(3,3,3)$ and $(b_1,b_2,b_3)=(1.5,1.5,1.5).$ The Bayes estimator of Lindley approximation under SEL function is $0.9032$ and under LINEX loss function is $0.9100$.  For MCMC, Bayes estimator under SEL function is $0.8813$ and under LINEX loss function is $0.8832$. Also, $95\%$ HPD interval of \rsk is $(0.68,0.99).$
\bibliography{MPRbib}
\end{document}